\documentclass[11pt]{amsart}

\usepackage[utf8]{inputenc}
\usepackage[english]{babel}
\usepackage{amsthm}

\usepackage{ amssymb,latexsym, amscd}
\usepackage[all]{xy}
\usepackage[square, numbers]{natbib}
\usepackage{graphicx}
\usepackage{pgf,tikz}
\usepackage{mathrsfs}
\usetikzlibrary{arrows}
\usepackage{mathrsfs}
\usepackage{amsmath}
\usepackage{color}
\usepackage{stmaryrd}
\usepackage{sseq}

\usepackage[T1]{fontenc}

\numberwithin{equation}{section}

\newtheorem{proposicion}{Proposition}[section]
\newtheorem{lema}[proposicion]{Lemma}
\newtheorem{teorema}[proposicion]{Theorem}
\newtheorem{corolario}[proposicion]{Corollary}

\theoremstyle{definition}
\newtheorem{observacion}[proposicion]{Remark}
\newtheorem{definicion}[proposicion]{Definition}
\newtheorem{ejemplo}[proposicion]{Example}

\newcommand{\bobs}{\begin{observacion}}
\newcommand{\eobs}{\end{observacion}}
\newcommand{\beq}{\begin{equation}}
\newcommand{\eeq}{\end{equation}}
\newcommand{\bprop}{\begin{proposicion}}
\newcommand{\eprop}{\end{proposicion}}
\newcommand{\blema}{\begin{lema}}
\newcommand{\elema}{\end{lema}}
\newcommand{\bejem}{\begin{ejemplo}}
\newcommand{\eejem}{\end{ejemplo}}
\newcommand{\bteo}{\begin{teorema}}
\newcommand{\eteo}{\end{teorema}}
\newcommand{\bdefin}{\begin{definicion}}
\newcommand{\edefin}{\end{definicion}}
\newcommand{\benum}{\begin{enumerate}}
\newcommand{\eenum}{\end{enumerate}}
\newcommand{\bcor}{\begin{corolario}}
\newcommand{\ecor}{\end{corolario}}
\newcommand{\bmat}{\begin{matrix}}
\newcommand{\emat}{\end{matrix}}
\newcommand{\barr}{\begin{array}}
\newcommand{\earr}{\end{array}}
\newcommand{\bcas}{\begin{cases}}
\newcommand{\ecas}{\end{cases}}
\newcommand{\bcen}{\begin{center}}
\newcommand{\ecen}{\end{center}}
\newcommand{\bdem}{\begin{proof}}
\newcommand{\edem}{\end{proof}}

\newcommand{\nin}{\noindent}
\newcommand{\tit}{\textit}

\newcommand{\tbf}{\textbf}

\newcommand{\mbb}{\mathbb}

\newcommand{\ol}{\overline}
\newcommand{\sse}{\subseteq}
\newcommand{\sm}{\setminus}
\newcommand{\vn}{\varnothing}

\newcommand{\PP}{{\mathcal P}}

\newcommand{\integer}{\ensuremath{{\mathbb Z}}}

\newcommand{\real}{\ensuremath{{\mathbb R}}}

\begin{document}
\title{The evasiveness conjecture and graphs on $2p$ vertices}
\author{Andr\'es Angel, Jerson Borja}
\email{ja.angel908@uniandes.edu.co,jm.borja583@uniandes.edu.co}

\maketitle
\begin{abstract} The Evasiveness conjecture have been proved for properties of graphs on a prime-power number of vertices and the six vertices case. The ten vertices case is still unsolved. In this paper we study the size of the automorphism group of a graph on $2p$ vertices to estimate the Euler characteristic of monotone non-evasive graph properties and get some conditions such graph properties must satisfy. We also do this by means of Oliver groups and give some lower bounds for the dimension of the simplicial complex associated to a nontrivial monotone non-evasive graph property. We apply our results to graphs on ten vertices to get conditions on potential counterexamples to the evasiveness conjecture in the ten vertices case.  
\end{abstract}

\section{Introduction}\label{eva-intro}

A \tit{graph property} $\PP$ is a collection of graphs on $n$ vertices that is closed under isomorphism of graphs. The complexity of $\PP$, denoted $c(\PP)$, is the minimum $k$ between 0 and ${n\choose 2}$ such that in the worst we are forced to ask $k$ questions of the form \tit{is $\{i,j\}$ an edge of $G$?} to an oracle that answers yes or no, in order to determine if $G$ belongs to $\PP$. When $c(\PP)={n\choose 2}$ we say that the property $\PP$ is \tit{evasive}, otherwise $\PP$ is \tit{non-evasive}.\\

We say a graph property $\PP$ is \tit{monotone} if it is closed under removal of edges. The property $\PP$ is called \tit{trivial} if it is either empty or is the family of all subsets of ${V\choose 2}$, otherwise $\PP$ is called \tit{nontrivial}.\\

The following conjecture was proposed by R. Karp. \\

\tit{The evasiveness conjecture for graph properties:} every nontrivial monotone graph property is evasive.\\

The ingenious approach of Kahn, Saks and Sturtevant \cite{kahn} permitted to solve the prime-power case and gave new directions on how to attack various related problems to the evasiveness conjecture.

In their approach, Kahn \tit{et al} associate to each nonempty monotone graph property $\PP$ a simplicial complex $\Delta\PP$ on the set of two-element subsets of $\{1,2,\ldots,n\}$. A collection of two-elements subset of $\{1,2,\ldots,n\}$ is a face of $\Delta \PP$ if and only if that collection corresponds to a graph belonging to $\PP$. By abuse of notation we denote both the graph property and the associated simplicial complex by $\PP$. 

The concept of evasiveness can be defined for simplicial complexes in a similar way to evasiveness of graph properties. In fact, a simplicial complex $K$ can be defined to be \tit{non-evasive} in an inductive way as follows. Define $K$ to be non-evasive if either $K$ is a single vertex or if there is a vertex $v$ in $K$ such that both the link $lk_K(v)$ and the deletion $del_K(v)$ are non-evasive (for a vertex $v$ in $K$, the link of $v$ is $lk_K(v)=\{A\sse V(K)\setminus\{v\}:\ A\cup\{v\}\in K\}$ and the deletion of $v$ is $del_K(v)=\{A\sse V\setminus\{v\}:\ A\in K\}$). With either definition of non-evasive simplicial complex we have 

\bteo\label{collapsible}(\tbf{Kahn-Saks-Sturtevant, \cite{kahn}}) A non-evasive complex is collapsible.
\eteo

The immediate consequence of this theorem is that a non-evasive complex is $\integer$-acyclic and therefore, $\integer/p$-acyclic for every prime number $p$. We also have that the Euler characteristic of an non-evasive simplicial complex is 1. 

With the help of results of R. Oliver \cite{oliver}, theorem \ref{collapsible} implies 

\bteo\label{primepower}(\tbf{Kahn-Saks-Sturtevant, \cite{kahn}}) The evasiveness conjecture is true for properties of graphs on a prime-power number of vertices. 
\eteo
 
In this paper we study the sizes of automorphism groups of graphs to estimate $\chi(\PP)$, where $\PP$ is a monotone property of graphs on $2p$ vertices where $p$ is prime. First, for a nonempty monotone property of graphs on $n$ vertices $\PP$ we can describe the Euler characteristic of $\PP$ in the form
$$\chi(\PP)=\sum_{[G]\sse \PP}(-1)^{m_G-1}|[G]|,$$
where $m_G$ is number of edges of $G$ and $|[G]|$ is the size of the isomorphism class $[G]$. Then, we describe all graphs on $p$ and $2p$ vertices $G$ for which $p$ does not divide $|[G]|$ (section \ref{pvertices}).\\
If we assume $\PP$ to be non-evasive, then $\chi(\PP)=1$ and this implies that $\PP$ must contain some of the graphs described in lemma \ref{D5} (below).

With the use of some Oliver groups and a result of Oliver \cite{oliver} (theorem \ref{th2} below), we show that $\PP$ necessarily has to contain some particular classes of graphs. We do this for graphs on $2p$ vertices and for graphs on $p^r+t$ vertices where $(p^r-1,t)=1$ (section \ref{eva-oliver} below). Using this results we can give lower bounds for the dimension of $\PP$ nontrivial monotone and non-evasive. In the particular case of a nontrivial monotone and non-evasive property $\PP$ of graphs on $2p$ vertices we get that $dim \PP\geq 4p-1$. This bound improves in the case of $2p$ vertices a general lower bound given by Bj\"orner for the dimension of a vertex homogeneous simplicial complex $K$ on a finite set of cardinality $m$ with $\chi(K)=1$ (see \cite{lutz}).

The applications of our results in the 10 vertices case gives us a reduction of the problem of evasiveness and leave us with conditions for potential counterexamples to the evasiveness conjecture in the 10 vertices case.  

\section{Simplicial complexes, Euler characteristic and automorphism group of graphs}

Let $V$ be a finite set. An (abstract) \tit{simplicial complex} on $V$ is a collection $K$ of subsets of $V$ such that (i) $\{v\}\in K$ for all $v\in V$ and $(ii)$ $A\in K$ and $B\sse A$ implies $B\in K$. If $A\in K$, $A$ is a \tit{face} of $K$. The dimension of $A$ is $|A|-1$. When the set of vertices $V$ is a face of $K$ we say that $K$ is a \tit{simplex}.\\

The \tit{geometric realization} of $K$ ,denoted $|K|$ can be construct as follows. If $V=\{v_1,v_2,\ldots,v_n\}$, identify $v_i$ with the standard basis vector $e_i\in\real^n$, then $|K|$ is obtained as the union of all convex hulls $\langle A\rangle=conv\{e_i:\ v_i\in A\}$ for $A\in K$.\\
If $K$ has $f_i$ faces of dimension $i$, then the Euler characteristic of $K$ is 
\[\chi(K)=\sum_{i\geq 0}(-1)^if_i.\]

The \tit{automorphism group} of $K$, denoted $Aut(K)$, is the collection of all permutations of $V$ which leave $K$ invariant. If $\Gamma$ is a subgroup of $Aut(K)$, then $\Gamma$ acts on $|K|$ by extending linearly the action on vertices. We write $|K|^ {\Gamma}$ for the fixed points of this action. The space $|K|^{\Gamma}$ is describe in an abstract way as follows. Define a simplicial complex $K^{\Gamma}$ whose vertices are the orbits of the action of $\Gamma$ on $V$ that are faces of $K$ and $\{A_1,A_2,\ldots,A_r\}$ is a face of $K^{\Gamma}$ if and only if $A_1\cup A_2\cup\cdots \cup A_r$ is a face of $K$.\\
If we identify each vertex $A_i$ of $K^{\Gamma}$ with the barycenter of $|A_i|$ in $|K|$, then the geometric realization of $K^{\Gamma}$ is just $|K|^{\Gamma}$.\\

If $\PP$ is a nonempty monotone and non-evasive graph property, then $\PP$ considered as a simplicial complex is collapsible and this implies that $\chi(\PP)=1$. Therefore, if $\chi(\PP)\neq 1$, then $\PP$ being monotone and nontrivial is evasive. This suggests the idea of examining the Euler characteristic of monotone graph properties of graphs on $n$ vertices.\\

For a given graph $G$ on $n$ vertices let $[G]$ denote its isomorphism class and for any pair of graphs $G,H$ on $n$ vertices let us write $[G]\leq [H]$ if and only if $G$ is isomorphic to some subgraph of $H$. We can write the Euler characteristic $\chi(\PP)$ as $$\chi(\PP)=\sum_{[G]\sse \PP}(-1)^{m_G-1}|[G]|,$$
where $\PP$ is supposed to be nonempty, $m_G$ represents the number of edges of the graph $G$ (so that $G$ corresponds to a face of dimension $m_G-1$ of the simplicial complex $\PP$), $|[G]|$ is the size of the isomorphism class $[G]$ and the sum is taken over all isomorphism classes of graphs contained in $\PP$. The idea is that in many instances there is a common divisor $d>1$ of all the sizes $|[G]|$ for $G\in\PP$. As a result, $d$ divides $\chi(\PP)$ so that $\chi(\PP)\neq 1$ and we conclude that $\PP$ is evasive.\\
Our observation is that if $\PP$ is nontrivial monotone and non-evasive, then $\PP$ must contain some graphs $G$ for which $d$ is not a divisor of $|[G]|$.\\

In this paper we are going to use the notation of Harary \cite{harary} for complements, unions and joins of graphs.

For a graph $G$ with vertex set $V(G)$ and edge set $E(G)$, its complement $\ol{G}$ has $V(\ol{G})=V(G)$ and $E(\ol{G})={V(G)\choose 2}\sm E(G)$. For graphs $G_1$ and $G_2$ with $V(G_1)\cap V(G_1)=\vn$, the graph $G_1\cup G_2$ has $V(G_1\cup G_2)=V(G_1)\cup V(G_2)$ and $E(G_1\cup G_2)=E(G_1)\cup E(G_2)$. If $V(G_1)\cap V(G_2)=\vn$, the \tit{join} graph $G_1+G_2$ has $V(G_1+ G_2)=V(G_1)\cup V(G_2)$ and $E(G_1+ G_2)=E(G_1)\cup E(G_2)\cup \{\{v_1,v_2\}:v_1\in V(G_1), v_2\in V(G_2)\}$. For instance, $\ol{K}_n$ is the graph on $n$ vertices without edges and $K_{n,m}=\ol{K}_m+\ol{K}_n$. If the graph $G$ is the union of $k$ connected components isomorphic to a graph $H$, then $G$ is denoted by $kH$.\\

In order to study the divisors of $|[G]|$, we observe that $|[G]|=\frac{n!}{|Aut(G)|}$, where $Aut(G)$ denotes the automorphism group of $G$. Therefore, we investigate the divisors of $|Aut(G)|$.\\ 

The following basic properties about automorphism groups can be found in \cite{harary}, chapter 14. 

\blema\label{desc} (i)If $G$ is a connected graph, then $Aut(kG)$ is isomorphic to the wreath product $Aut(G)\wr S_k$.\\
(ii) If $G_1$ and $G_2$ are disjoint and connected non-isomorphic graphs, then $Aut(G_1\cup G_2)\cong Aut(G_1)\times Aut(G_2)$.\\
(iii) If no component of $\ol{G}_1$ is isomorphic to a component of $\ol{G}_2$, then
\[Aut(G_1+G_2)\cong Aut(G_1)\times Aut(G_2).\]
(iv) If $G$ is a graph on $n$ vertices, then the group $Aut(G)$ decomposes as
\[
Aut(G)\cong (Aut(G_1)\wr S_{n_1})\times\cdots\times (Aut(G_s)\wr S_{n_s})
\]
where the $G_i$'s are the distinct connected components of $G$ and $n_i$ is the number of components of $G$ isomorphic to $G_i$. If $m_i$ is the number of vertices of $G_i$, then $n=n_1m_1+n_2m_2+\cdots +n_sm_s$, $Aut(G_i)$ is isomorphic to a subgroup of $S_{m_i}$ and $|Aut(G)|$ divides $\prod_in_i!\cdot(m_i!)^{n_i}$.
\elema

\blema \label{D1}
Suppose that a graph $G$ on $n$ vertices has exactly $k$ vertices of a fixed degree $r$, where $0\leq k\leq n$. Then $Aut(G)$ is isomorphic to a subgroup of $S_k\times S_{n-k}$. In particular, ${n\choose k}$ divides $|[G]|$.
\elema

This is because every element of $Aut(G)$ preserves the set of vertices of degree $r$ and also preserves the set of vertices that are not of degree $r$. Then $|Aut(G)|$ divides $k!(n-k)!$ and $\frac{n!}{k!(n-k)!}$ divides $|[G]|$.\\

A graph $G$ is called \tit{regular} if all of its vertices have the same degree. If such degree is $r$, we say that $G$ is $r$-regular. 
When studying the divisors of the size of automorphism groups of regular graphs, the following result of N. Wormald is very useful. 

\bteo\label{Wormald} \tbf{(Wormald, \cite{wormald})} Let $G$ be a connected $r$-regular graph on $n$ vertices, where  $r>0$. Then $|Aut(G)|$ divides 
\[
rn\prod_pp^{\beta}
\]
where the product is taken over all prime numbers $p\leq r-1$ and $\beta$ is given by  
\[
\sum_{p^{\alpha}\leq r-1}\left\lfloor\frac{n-2}{p^{\alpha}}\right\rfloor
\]
\eteo

\bcor \label{regular3} In the hypothesis of theorem \ref{Wormald}, if $r<3$, then $|Aut(G)|$ divides $rn$.
\ecor

\section{Graphs on $p$ and $2p$ vertices}\label{pvertices}

In this section we describe all graphs $G$ on $p$ vertices and also on $2p$ vertices for which $p$ is not a divisor of $|[G]|$. We assume that $p$ is an odd prime. First, we deal with graphs on $p$ vertices.\\

In the hypothesis of lemma \ref{D1} we have that ${p\choose k}$ divides $|[G]|$. Since $p$ is prime, $p$ divides ${p\choose k}$ for $0<k<p$ so the only cases in which $p$ may not be a divisor of $|[G]|$ are $k=p$ ($G$ is a regular graph) and $k=0$ (there are no vertices of degree $r$ in $G$). When we vary the degree $r$, we conclude that if $p$ does not divide $|[G]|$, then $G$ is regular.  

\blema
If $G$ is a graph on $p$ vertices and $p$ does not divide $|[G]|$, then $G$ is a regular graph.
\elema

Since $|[G]|=\frac{p!}{|Aut(G)|}$, we have that $p$ does not divide $|[G]|$ if and only if $p$ divides $|Aut(G)|$, which means that $Aut(G)$ contains an element of order $p$. The elements of order $p$ in the symmetric group $S_p$ are precisely the $p$-cycles. We assume, without lost of generality, that the $p$-cycle $\sigma=(01\cdots p-1)$ is an automorphism of $G$, where the vertices $0,1,\ldots,p-1$ are regarded as the finite field $\mbb{F}_p$.\\

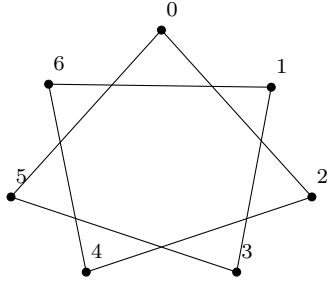
\begin{figure}
\bcen
\begin{tikzpicture}
\draw (0,4.22)-- (-2,2);
\draw (-2,2)-- (1,1);
\draw (1,1)-- (1.46,3.46);
\draw (1.46,3.46)-- (-1.5,3.5);
\draw (-1.5,3.5)-- (-1,1);
\draw (-1,1)-- (2,2);
\draw (2,2)-- (0,4.22);
\begin{scriptsize}
\draw [fill=black] (-1,1) circle (1.5pt);
\draw[color=black] (-0.86,1.28) node {$4$};
\draw [fill=black] (1,1) circle (1.5pt);
\draw[color=black] (1.14,1.28) node {$3$};
\draw [fill=black] (2,2) circle (1.5pt);
\draw[color=black] (2.14,2.28) node {$2$};
\draw [fill=black] (-2,2) circle (1.5pt);
\draw[color=black] (-1.86,2.28) node {$5$};
\draw [fill=black] (-1.5,3.5) circle (1.5pt);
\draw[color=black] (-1.36,3.78) node {$6$};
\draw [fill=black] (1.46,3.46) circle (1.5pt);
\draw[color=black] (1.6,3.74) node {$1$};
\draw [fill=black] (0,4.22) circle (1.5pt);
\draw[color=black] (0.14,4.5) node {$0$};
\end{scriptsize}
\end{tikzpicture}
\ecen
\caption{$C(2)$, $p=7$, $\sigma=(0123456)$.}
\label{hepta}
\end{figure}

Let us suppose that $G$ is not the graph $\ol{K}_p$. If the edge $\{r,s\}$ is in $G$ where $r<s$, reorder the $p$-cycle if necessary to assume that $r=0$. Thus, we can assume that the edge $\{0,s\}$ is in $G$. The $(s-1)$-th power of the $p$-cycle $\sigma$, $\sigma^{s-1}=(0\ s\ 2s\cdots (p-1)s)$, is again a $p$-cycle and an automorphism of $G$. We find that $G$ contains the edges $\{0,s\},\{s,2s\},\{2s,3s\},\ldots, \{(p-1)s,0\}$. These set of edges form a cycle graph of length $p$ that is a subgraph of $G$ and that we denote $C(s)$. With this notation we have $C(s)=C((p-1)s)$. Remember that the vertices $0,1,\ldots,p-1$ are the elements of $\mbb{F}_p$ and since $(p-1)s=-s=p-s$, we get $C(s)=C(p-s)$. In figure \ref{hepta} we illustrate the example for $p=7$, $s=2$, $\sigma=(0123456)$, $C(2)=C(5)$.

Every element $s\in\{1,2,\ldots,(p-1)/2\}$ such that $\{0,s\}\in G$ determines a cycle graph $C(s)$. Now, suppose that $\{0,s\}, \{0,r\}$ are in $G$, where $s,r\in\{1,2,\ldots,(p-1)/2\}$ and $r\neq s$. The cycles $C(s)$ and $C(r)$ have sets of edges $\{0,s\},\{s,{2s}\},\ldots,\{{(p-1)s},0\}$ and $\{0,r\},\{r,{2r}\},\ldots,\{{(p-1)r},0\}$ respectively.\\
We claim that $C(s)$ and $C(r)$ have no common edges. In fact, if we have $\{is,(i+1)s\}=\{jr,(j+1)r\}$, then there are two options: (i) $is=jr$, $is+s=jr+r$ which implies that $r=s$ and (ii) $is=jr+r$, $is +s=jr$, which implies $s=-r=k-r$ and this contradicts $s\in\{1,2,\ldots,(p-1)/2\}$.\\
Therefore, the graph $G$ can be decomposed as a union of disjoint cycle graphs of length $p$
\[
G=C(s_1)\cup C(s_2)\cup \cdots \cup C(s_l)
\]
where $\{0,{s_i}\}\in G$ and $s_i\in\{1,2,\ldots,(p-1)/2\}$. As a consequence, $G$ is an $l$-regular graph which has $lp$ edges.\\
We will denote the graph $C(s_1)\cup C(s_2)\cup \cdots \cup C(s_l)$ by $C(s_1,s_2\ldots,s_l)$. The complement graph of a graph of the form $C(s_1,s_2\ldots,s_l)$ also has the form $C(s_1,s_2\ldots,s_l)$ or is equal to the graph $\ol{K_p}$. All the cycle graphs $C(s)$ are isomorphic, in fact, given $s\in\{1,2,\ldots,(p-1)/2\}$, the map $\mbb{F}_p\to\mbb{F}_p$ given by $x\mapsto sx$ gives an isomorphism between $C(1)$ and $C(s)$.\\

Each cycle graph $C(s_i)$ is fixed under the action of the $p$-cycle $\sigma=(01\cdots {p-1})$ and as a consequence $G$ remains fixed under the action of $\sigma$, that is, $\sigma\in Aut(G)$. We have proved the following.

\blema\label{D3} If $G$ is a graph on $p$ vertices, then $p$ does not divide $|[G]|$ if and only if $G=\ol{K_p}$ or $G$ is isomorphic to one of the graphs $C(s_1,s_2\ldots,s_l)$.
\elema

Now we go to graphs on $2p$ vertices. If $G$ is a graph on $2p$ vertices ($p$ an odd prime) and $G$ have exactly $k$ vertices of a fixed degree $r$, then ${2p\choose k}$ divides $|[G]|$ by lemma \ref{D1}. The only way $p$ does not divide ${2p\choose k}=\frac{(2p)!}{k!(2p-k)!}$ is that $p^2$ divides $k!(2p-k)!$, which occurs if and only if $k=p$ or $k=2p$. Thus, we assume that $p$ does divide $|[G]|$ and consider the following two cases.\\

\nin\tbf{Case 1.} $k=p$. In this case there are exactly $p$ vertices of degree $r$, say they are $a_1,\ldots,a_p$. The other $p$ vertices, say $b_1,\ldots,b_p$, must have the same degree $s$, where $s\neq r$. Each automorphism of $G$ preserves the $a_i's$ and also preserves the $b_j's$. The group $Aut(G)$ is isomorphic to a subgroup of $S_{\{a_i\}}\times S_{\{b_j\}}\cong S_p\times S_p$.\\
Let $G_1$ be the subgraph of $G$ which is the subgraph \tit{induced} by the vertices $a_1,\ldots,a_p$, $G_2$ be the subgraph of $G$ induced by the vertices $b_j$'s and let $H$ be the subgraph of $G$ whose edges are all edges of $G$ having the form $\{a_i, b_j\}$ (see figure \ref{GG1G2H}). 

\begin{figure}
\bcen
\begin{tikzpicture}

\draw (-0.25,-.75)-- (-.75,-.5);   
\draw (-0.25,-.75)-- (-0.25,.75);  
\draw (-.75,-.5)-- (-1,0);           
\draw (-1,0)-- (-.75,.5);            
\draw (-.75,.5)-- (-0.25,.75);         

\draw (0.25,.75)-- (1,0);        
\draw (0.25,.75)-- (.75,-.5);   
\draw (.75,.5)-- (.75,-.5);       
\draw (.75,.5)-- (0.25,-.75);    
\draw (1,0)-- (0.25,-.75);        

\draw (-0.25,-.75)-- (-1,0);        
\draw (-0.25,-.75)-- (-.75,.5);     
\draw (-.75,-.5)-- (-.75,.5);        
\draw (-.75,-.5)-- (-0.25,.75);    
\draw (-1,0)-- (-0.25,.75);         

\draw (0.25,.75)-- (-0.25,-.75);
\draw (0.25,.75)-- (-.75,-.5);   
\draw (0.25,.75)-- (-1,0);       
\draw (0.25,.75)-- (-.75,.5);    
\draw (0.25,.75)-- (-0.25,.75);   
\draw (.75,.5)-- (-0.25,-.75);   
\draw (.75,.5)-- (-.75,-.5);      
\draw (.75,.5)-- (-1,0);           
\draw (.75,.5)-- (-.75,.5);       
\draw (.75,.5)-- (-0.25,.75);    
\draw (1,0)-- (-0.25,-.75);       
\draw (1,0)-- (-.75,-.5);          
\draw (1,0)-- (-1,0);              
\draw (1,0)-- (-.75,.5);           
\draw (1,0)-- (-0.25,.75);        
\draw (.75,-.5)-- (-0.25,-.75);   
\draw (.75,-.5)-- (-.75,-.5);      
\draw (.75,-.5)-- (-1,0);           
\draw (.75,-.5)-- (-.75,.5);        
\draw (.75,-.5)-- (-0.25,.75);     
\draw (0.25,-.75)-- (-0.25,-.75); 
\draw (0.25,-.75)-- (-.75,-.5);    
\draw (0.25,-.75)-- (-1,0);         
\draw (0.25,-.75)-- (-.75,.5);      
\draw (0.25,-.75)-- (-0.25,.75);   

\begin{scriptsize}
\draw [fill=black] (0.25,.75) circle (1pt);
\draw (0.287,.98) node {$b_1$};
\draw [fill=black] (.75,.5) circle (1pt);
\draw (1.,.617) node {$b_2$};
\draw [fill=black] (1,0) circle (1pt);
\draw (1.215,0) node {$b_3$};
\draw [fill=black] (.75,-.5) circle (1pt);
\draw (1.,-.617) node {$b_4$};
\draw [fill=black] (0.257,-.75) circle (1pt);
\draw (0.287,-.98) node {$b_5$};
\draw [fill=black] (-0.25,-.75) circle (1pt);
\draw (-0.287,-.98) node {$a_5$};
\draw [fill=black] (-.75,-.5) circle (1pt);
\draw (-1,-.617) node {$a_4$};
\draw [fill=black] (-1,0) circle (1pt);
\draw (-1.215,0) node {$a_3$};
\draw [fill=black] (-.75,.5) circle (1pt);
\draw (-1,.617) node {$a_2$};
\draw [fill=black] (-0.25,.75) circle (1pt);
\draw (-0.287,.98) node {$a_1$};
\draw (0,-1.4) node {$G$};
\end{scriptsize}
\end{tikzpicture}
\ecen

\bcen
\begin{tikzpicture}

\draw (-0.25,-.75)-- (-.75,-.5);   
\draw (-0.25,-.75)-- (-0.25,.75);  
\draw (-.75,-.5)-- (-1,0);           
\draw (-1,0)-- (-.75,.5);            
\draw (-.75,.5)-- (-0.25,.75);         

\draw (-0.25,-.75)-- (-1,0);        
\draw (-0.25,-.75)-- (-.75,.5);     
\draw (-.75,-.5)-- (-.75,.5);        
\draw (-.75,-.5)-- (-0.25,.75);    
\draw (-1,0)-- (-0.25,.75);         

\begin{scriptsize}
\draw [fill=black] (-0.25,-.75) circle (1pt);
\draw (-0.287,-.98) node {$a_5$};
\draw [fill=black] (-.75,-.5) circle (1pt);
\draw (-1,-.617) node {$a_4$};
\draw [fill=black] (-1,0) circle (1pt);
\draw (-1.215,0) node {$a_3$};
\draw [fill=black] (-.75,.5) circle (1pt);
\draw (-1,.617) node {$a_2$};
\draw [fill=black] (-0.25,.75) circle (1pt);
\draw (-0.287,.98) node {$a_1$};
\draw (-0.6,-1.4) node {$G_1$};
\end{scriptsize}
\end{tikzpicture}
\begin{tikzpicture}

\draw (0.25,.75)-- (1,0);        
\draw (0.25,.75)-- (.75,-.5);   
\draw (.75,.5)-- (.75,-.5);       
\draw (.75,.5)-- (0.25,-.75);    
\draw (1,0)-- (0.25,-.75);        

\begin{scriptsize}
\draw [fill=black] (0.25,.75) circle (1pt);
\draw (0.287,.98) node {$b_1$};
\draw [fill=black] (.75,.5) circle (1pt);
\draw (1.,.617) node {$b_2$};
\draw [fill=black] (1,0) circle (1pt);
\draw (1.215,0) node {$b_3$};
\draw [fill=black] (.75,-.5) circle (1pt);
\draw (1.,-.617) node {$b_4$};
\draw [fill=black] (0.257,-.75) circle (1pt);
\draw (0.287,-.98) node {$b_5$};
\draw (0.6,-1.4) node {$G_2$};
\end{scriptsize}
\end{tikzpicture}
\begin{tikzpicture}

\draw (0.25,.75)-- (-0.25,-.75);
\draw (0.25,.75)-- (-.75,-.5);   
\draw (0.25,.75)-- (-1,0);       
\draw (0.25,.75)-- (-.75,.5);    
\draw (0.25,.75)-- (-0.25,.75);   
\draw (.75,.5)-- (-0.25,-.75);   
\draw (.75,.5)-- (-.75,-.5);      
\draw (.75,.5)-- (-1,0);           
\draw (.75,.5)-- (-.75,.5);       
\draw (.75,.5)-- (-0.25,.75);    
\draw (1,0)-- (-0.25,-.75);       
\draw (1,0)-- (-.75,-.5);          
\draw (1,0)-- (-1,0);              
\draw (1,0)-- (-.75,.5);           
\draw (1,0)-- (-0.25,.75);        
\draw (.75,-.5)-- (-0.25,-.75);   
\draw (.75,-.5)-- (-.75,-.5);      
\draw (.75,-.5)-- (-1,0);           
\draw (.75,-.5)-- (-.75,.5);        
\draw (.75,-.5)-- (-0.25,.75);     
\draw (0.25,-.75)-- (-0.25,-.75); 
\draw (0.25,-.75)-- (-.75,-.5);    
\draw (0.25,-.75)-- (-1,0);         
\draw (0.25,-.75)-- (-.75,.5);      
\draw (0.25,-.75)-- (-0.25,.75);   

\begin{scriptsize}
\draw [fill=black] (0.25,.75) circle (1pt);
\draw (0.287,.98) node {$b_1$};
\draw [fill=black] (.75,.5) circle (1pt);
\draw (1.,.617) node {$b_2$};
\draw [fill=black] (1,0) circle (1pt);
\draw (1.215,0) node {$b_3$};
\draw [fill=black] (.75,-.5) circle (1pt);
\draw (1.,-.617) node {$b_4$};
\draw [fill=black] (0.257,-.75) circle (1pt);
\draw (0.287,-.98) node {$b_5$};
\draw [fill=black] (-0.25,-.75) circle (1pt);
\draw (-0.287,-.98) node {$a_5$};
\draw [fill=black] (-.75,-.5) circle (1pt);
\draw (-1,-.617) node {$a_4$};
\draw [fill=black] (-1,0) circle (1pt);
\draw (-1.215,0) node {$a_3$};
\draw [fill=black] (-.75,.5) circle (1pt);
\draw (-1,.617) node {$a_2$};
\draw [fill=black] (-0.25,.75) circle (1pt);
\draw (-0.287,.98) node {$a_1$};
\draw (0,-1.4) node {$H$};
\end{scriptsize}
\end{tikzpicture}
\ecen

\caption{}
\label{GG1G2H}
\end{figure}

If $p^2$ divides the size of $Aut(G)$, then $Aut(G)$ contains a subgroup of order $p^2$. The subgroups of order $p^2$ in $S_p\times S_p$ have the form $\langle\sigma\rangle\times\langle\tau\rangle$  where both $\sigma$ and $\tau$ have order $p$.  Thus, $Aut(G)\sse S_{\{a_i\}}\times S_{\{b_j\}}$ contains a subgroup of the form $\langle\alpha\rangle\times\langle\beta\rangle$ where both elements $\alpha$ and $\beta$ have order $p$ ($\alpha$ and $\beta$ are $p$-cycles).\\
Let $t$ be the number of the $b_j's$ that are adjacent to $a_1$. Then, for each $i=1,2,\ldots,p$, the number of the $b_j's$ adjacent to $a_i$ is also $t$ because automorphisms of $G$ preserve the adjacency and $\alpha\times id$ acts transitively on the $a_i's$ while $\alpha\times id$ fixes the $b_j$'s. Thus, the number of edges of $H$ is $m_H=pt$. Similarly, if $u$ is the number of $a_i's$ adjacent to $b_1$, then $m_H=pu$. Hence, $u=t$.\\

If $t=0$, then $H$ has no edges. In this case $G=G_1\cup G_2$ and since no component of $G_1$ is isomorphic to a component of $G_2$, we have that $Aut(G)\cong Aut(G_1)\times Aut(G_2)$. Therefore, $p$ does not divide $|[G]|$ if and only if $p^2$ divides $|Aut(G)|$ and this is equivalent to say that $p$ divides both $|Aut(G_1)|$ and $|Aut(G_2)|$. Both graphs $G_1$ and $G_2$ have $p$ vertices and $p$ divides the size of their automorphism groups, so both graphs $G_1$ and $G_2$ have the form of the graphs described in lemma \ref{D3}. Then, $G_1$ and $G_2$ are regular graphs of different degree.\\

If $t>0$, the permutation $id\times \beta$ permutes transitively the $b_j$'s and fixes $a_1$, this implies that $t=p$, $H=K_{p,p}$ and $G=G_1+G_2$. Therefore, $\ol{G}=\ol{G_1}\cup\ol{G_2}$ that lies in the case $t=0$ above, so that $\ol{G}_1$ and $\ol{G}_2$ and also $G_1$ and $G_2$ have the form of of the graphs described in lemma \ref{D3}.\\

\nin\tbf{Case 2.} $k=2p$. In this case $G$ is a regular graph of degree $r$. When $r=0$, $G=\ol{K}_{2p}$. We assume $r>0$. By lemma \ref{desc}, the group $Aut(G)$ decomposes as
\[
Aut(G)\cong (Aut(G_1)\wr S_{n_1})\times\cdots\times (Aut(G_s)\wr S_{n_s}), 
\]
where the $G_i$'s are the distinct connected components of $G$ and $n_i$ is the number of components of $G$ isomorphic to $G_i$. If  $m_i$ is the number of vertices of $G_i$, then $2p=n_1m_1+n_2m_2+\cdots +n_sm_s$, $Aut(G_i)$ is isomorphic to a subgroup of $S_{m_i}$ and $|Aut(G)|$ divides $\prod_in_i!\cdot(m_i!)^{n_i}$. Since there are no isolated vertices in $G$, we have that $m_i\geq 2$ for all $i$.

We want to determine conditions under which $p^2$ divides $\prod_in_i!\cdot(m_i!)^{n_i}$.  If some $n_i\geq p$, then $2p\geq n_im_i\geq 2n_i\geq 2p$ which implies $s=1$, $n_1=p$ and $m_1=2$. It follows that $G$ is the union of $p$ disjoint copies of the complete graph on 2 vertices $K_2$, besides $|Aut(G)|=p!2^p$ is not divisible by $p^2$.\\

If $n_i<p$ for all $i$, then the only cases in which $p^2$ divides $\prod_in_i!\cdot(m_i!)^{n_i}$ are the following:\\
\nin (i) $m_i=2p$ for some $i$. In this case $G$ is a regular connected graph of degree $r>0$.\\
\nin (ii) $m_i=p$ and $m_j=p$ for $i\neq j$. In this case  $G$ is the union of two non-isomorphic regular connected graphs on $p$ vertices, both of degree $r$.\\	
\nin (iii) $m_i=p$ and $n_i= 2$ for some $i$. $G$ is the union of two copies of a regular connected graph on $p$ vertices of degree $r$.\\
In the case that $G$ is regular and connected of degree $r>0$ (case (i)), consider $\ol{G}$ which is also regular of degree $2p-r-1$ and $|Aut(G)|=|Aut(\ol{G})|$. If $\ol{G}$ is also connected, take the one between $G$ and $\ol{G}$ that has fewer degree, so we can assume that $G$ is connected of degree $r\leq \frac{2p-1}{2}$. This condition on $r$ is equivalent to $r<p$. Then we apply Wormald's theorem (theorem \ref{Wormald}) to obtain that $|Aut(G)|$ divides $r(2p)\prod_{q\leq r-1}q^{\beta}$ and we see from this that $p^2$ is not a divisor of $|Aut(G)|$ (although $p$ could be a divisor of $|Aut(G)|$).\\
On the other hand, if $\ol{G}$ is not connected and $p^2$ divides $|Aut(G)|=|Aut(\ol{G})|$, then $\ol{G}$ lies in case (ii) or case (iii). Thus, $\ol{G}$ is the union of two regular graphs on $p$ vertices of the same degree, $\ol{G}=H_1\cup H_2$ and $G=\ol{H}_1+\ol{H}_2$.\\

We conclude that the only way that $p^2$ divides $|Aut(G)|$, where $G$ is a regular graph of degree $r$ is that $G$ is of the form $G_1\cup G_2$ or $G_1+G_2$, where $G_1$ and $G_2$ have the form of the graphs described in lemma \ref{D3}.\\

We resume our classification in the following lemma:

\blema\label{D5}
 Let $G$ be a graph on $2p$ vertices such that $p^2$ divides $|Aut(G)|$, or equivalently that $p$ does not divide $|[G]|$. Then, $G$ is isomorphic to a graph of the form $G_1\cup G_2$ or to a graph of the form $G_1+G_2$, where $G_1$ and $G_2$ are graphs on $p$ vertices of the form described in lemma \ref{D3}. The graphs $G_1$ and $G_2$ are allowed to be isomorphic.
\elema

\bobs
To determine all graphs $G$ on $2p$ vertices such that $p$ does not divide $|[G]|$, it suffices to determine all graphs $H$ on $p$ vertices such that $p$ does not divide $|[H]|$ (lemma \ref{D3} contains the precise description of such graphs $H$) and then take unions $G_1\cup G_2$ and joins $G_1+G_2$ of such graphs. Graphs of the form $G_1+G_2$ are obtained from those of the form $G_1\cup G_2$ by the relation $\ol{G_1\cup G_2}=\ol{G_1}+\ol{G_2}$.
\eobs

\section{Oliver groups and lower bounds for the dimension of non-evasive monotone graph properties}\label{eva-oliver}

A finite group $\Gamma$ is an \tit{Oliver} group if $\Gamma$ has a normal $p$-subgroup $\Gamma_1$ such that the quotient $\Gamma/\Gamma_1$ is cyclic. Thus, all finite $p$-groups are oliver groups. The importance of these groups is justified by the following theorem.
 
\bteo \label{th2}(\tbf{Oliver}) Let $K$ be a simplicial complex, $\Gamma$ be a finite subgroup of $Aut(K)$ and $p$ be a fixed prime. Assume that\\
i) $|K|$ is $\integer/p$-acyclic and\\
ii) $\Gamma$ has a normal $p$-subgroup $\Gamma_1$ such that the quotient $\Gamma/\Gamma_1$ is cyclic.\\
Then $\chi(|K|^{\Gamma})=1$. In particular, $|K|^{\Gamma}$ is nonempty. 
\eteo

Let $\PP$ be a nonempty monotone and non-evasive graph property. If $\Gamma$ is an Oliver group acting on $|\PP|$, where $\PP$ is regarded as simplicial complex, then $|\PP|^{\Gamma}\neq \vn$. In the abstract version of $|\PP|^{\Gamma}$ we obtain that $\PP^{\Gamma}\neq\vn$ and this means that there are some graphs in $\PP$ that are also (union of) orbits of the action of $\Gamma$ on the set of vertices of $\PP$ (the 2-elements subsets of $\{1,2,\ldots,n\}$).

We will consider Oliver groups that are subgroups of $S_n$ because they automatically act on $\PP$ as automorphisms of $\PP$. We conclude that at least one of the orbits of $\Gamma$ acting on the 2-element subsets of $\{1,2,\ldots,n\}$ must be a graph belonging to $\PP$.\\

For the first example, suppose $n=p^r$ is a prime power and consider $GF(p^r)$, the finite field of $p^r$ elements. For $a,b\ GF(p^r), a\neq 0$ let $\phi_{a,b}:GF(p^r)\to GF(p^r)$ be the affine linear map defined by $\phi_{a,b}(x)=ax+b$. If $\Gamma_{p^r}=\{\phi_{a,b}:\ a,b\ GF(p^r), a\neq 0\}$, then $\Gamma_{p^r}$ is an Oliver group. In fact, this group $\Gamma_{p^r}$ is the one used by Kahn, Saks and Sturtevant in the proof of theorem \ref{primepower} (see \cite{kahn}). The normal $p$-subgroup with cyclic quotient is $\Gamma_1=\{\phi_{1,b}: b\in GF(p^r)\}$. The group $\Gamma_{p^r}$ acts doubly transitively on $GF(p^r)$ (the set of vertices) which implies that there is just one orbit of $\Gamma_{p^r}$ acting on the two-element subsets of $\{1,2,\ldots, p^r\}$, the complete graph $K_{p^r}$.\\

Now suppose that $n=p^r+t$ and consider the group $\Gamma_{p^r}\times \integer/t$ where $\Gamma_{p^r}$ acts on the vertices $1,2,\ldots,p^r$ as above while fixes the vertices $p^r+1,\ldots,n$ and the factor $\integer/t$ acts on $1,2,\ldots,p^r$ trivially and acts on the remaining $t$ vertices $p^r+1,\ldots,n$ by permuting them cyclically.

We know that $\Gamma_{p^r}$ has the normal $p$-subgroup $\Gamma_1$ with cyclic quotient. Then $\Gamma_1$ is also a normal $p$-subgroup of $\Gamma\times \integer/t$ with quotient isomorphic to $(\Gamma/\Gamma_1)\times \integer/t$ which is cyclic when $p^r-1$ and $t$ are coprimes. Thus, the group $\Gamma\times\integer/t$ is an Oliver group if $p^r-1$ and $t$ are coprimes.

The orbits of the action of $\Gamma_{p^r}\times \integer/t$ on the two element subsets of $\{1,2,\ldots,n\}$ are the complete graph on the vertices $1,2,\ldots,p^r$ with the remaining $t$ vertices isolated, that is, $K_{p^r}\cup \ol{K}_t$; the complete bipartite graph $K_{p^r,t}$ and the orbits of the form $\ol{K}_{p^r}\cup G$, where $G$ is one of the graphs on $t$ vertices $p^r+1,\ldots,n$ that are fixed under the action of $\integer/t$.

\bprop\label{oliverprop1}
Any nonempty monotone and non-evasive graph property on $n=p^r+t$ vertices, where $(p^r-1,t)=1$, has to contain some of the following graphs: $K_{p^r}\cup \ol{K}_{t}, K_{p^r,t}$ or one of the graphs of the form $\ol{K}_{p^r}\cup G$, where $G$ is one of the graphs on $t$ vertices $p^r+1,\ldots,n$ that are fixed under the action of $\integer/t$. 
\eprop

\bcor 
Any nonempty monotone and non-evasive property $\PP$ of graphs on $n=p^r+1$ vertices, where $n>3$, has to contain one of the two graphs $K_{p^r}\cup K_1$ or $K_{p^t,1}$. In particular $dim\ \PP\geq n-2$. Moreover, if $\PP$ is nontrivial then $\PP$ cannot contain both $K_{p^r}\cup K_1$ and $K_{p^t,1}$.
\ecor

\bdem
Applying proposition \ref{oliverprop1}, $\PP$ has to contain $K_{p^r}\cup K_1$ or $K_{p^t,1}$. These two graphs have $p^r(p^r-1)/2=(n-1)(n-2)/2$ and $p^r=n-1$ edges respectively. In any case $dim\ \PP\geq n-2$ if $n> 3$. The two graphs $A=K_{p^r}\cup K_1$ and $B=K_{p^t,1}$ are the only two orbits of the Oliver group $\Gamma\times \{\ast\}$ acting on the 2-element subsets of $1,2,\ldots,n$. Thus, $A$ and $B$ represent the only two (potential) vertices of the simplicial complex $\PP^{\Gamma\times \{\ast\}}$. If both $A$ and $B$ belong to $\PP$, then $A$ and $B$ are vertices of $\PP^{\Gamma\times \{\ast\}}$. By theorem \ref{th2}, $\chi(|\PP^{\Gamma\times \{\ast\}}|)=1$, which obligates $\PP^{\Gamma\times \{\ast\}}$ to contain the simplex $\{A;B\}$, which means that $K_{n}=A\cup B$ belongs to $\PP$ and this implies that $\PP$ is trivial. 
\edem

The following to results are generalizations of the corresponding results for 6 vertices in \cite{kahn}. In the proof of each result an Oliver group is used in combination with theorem \ref{th2}. Both results concern monotone non-evasive properties of graphs on $2p$ vertices.  

\bprop\label{matchings2p} Let $\PP$ be a nontrivial monotone and non-evasive graph property on $2p$ vertices, where $p$ is prime. Then all perfect matchings belong to $\PP$.
\eprop 

\bdem Let $\Gamma$ be the group generated by the permutations $(1\ p+1),(2\ p+2),\ldots, (p\ 2p) $ and $\alpha=(1\ 2\cdots\ p)(p+1\ p+2\cdots 2p)$. The subgroup $H$ of $\Gamma$ generated by $(1\ p+1),(2\ p+2),\ldots, (p\ 2p)$ is a normal 2-subgroup with quotient isomorphic to the subgroup of $\Gamma$ generated by $\alpha$, which is cyclic. Then $\Gamma$ is an Oliver group and $\PP^{\Gamma}\neq \vn$.

We claim that every orbit of $\Gamma$ acting on the two-element subsets of $\{1,2,\ldots,2p\}$ contains a perfect matching. One of such orbits is the set of edges $\{1,p\},\{2,p+1\},\ldots,\{p,2p\}$ (a perfect matching!).\\
Let $G$ be any other orbit and $\{r,s\}\in G$. We can assume that $r,s\in \{1,2,\ldots,p\}$, for if $r\leq p$ and $s>p$ then $s-p\neq r$ and the permutation $(s-p\ s)$ sends $\{r,s\}$ to $\{r,s-p\}$, and $s-p\in \{1,2,\ldots,p\}$. If $r,s>p$, then $(r-p\ r)(s-p\ s)$ sends $\{r,s\}$ to $\{r-p,s-p\}$ and both $r-p,s-p\in \{1,2,\ldots,p\}$. The orbit of $\{r,s\}$ under $\alpha$ is a cycle graph of length $p$ with vertices $\{1,2,\ldots,p\}$ which is a subgraph of $G$. Since $\alpha^{1-r}$ sends $r$ to 1, $G$ contains an edge of the form $\{1,t\}$ with $1<t\leq p$. Now, $(1\ p)$ sends $\{1,t\}$ to $\{p,t\}$ and so $\{p,t\}\in G$, then by considering the action of $\alpha$ on $\{p,t\}$, we conclude that $\{p,t\},\{p+1,t+1\},\{p+2,t+2\},\ldots, \{2p,t+p\}\in G$ (if necessary, when $t+i>p$ we subtract $p$ to obtain a value between 1 and $p$). This set of edges is a perfect matching and this ends the proof of our claim.

Since $\PP^{\Gamma}\neq \vn$, at least one of the orbits of $\Gamma$ belongs to $\PP$ and since $\PP$ is monotone, $\PP$ contains a perfect matching. Then, $\PP$ contains all perfect matchings for $\PP$ is closed under isomorphism of graphs.\edem

\bprop \label{2cpokpp} Let $\PP$ be a nontrivial monotone and non-evasive graph property on $2p$ vertices, where $p$ is and odd prime. Then at least one of $2C_p$, $K_{p,p}$ belongs to $\PP$.
\eprop

\bdem
The set of $2p$ vertices is going to be the union of two disjoint copies of the finite field $\mbb{F}_p=\{0,1,\ldots,p-1\}$, the second copy of $\mbb{F}_p$ will be labeled $\mbb{F}_p'=\{0',1',\ldots,(p-1)'\}$. Let $\Gamma$ be the group generated by the permutations $\alpha=(0\ 0')(1\ 1')\cdots (p-1\ (p-1)'), \beta=(01\cdots p-1), \gamma=(0'1'\cdots (p-1)')$. The subgroup of $\Gamma$ generated by $\beta$ and $\gamma$ is a normal $p$-subgroup of $\Gamma$ whose quotient is cyclic isomorphic to $\langle \alpha\rangle$. Thus, $\Gamma$ is an Oliver group and $\PP^{\Gamma}\neq \vn$.

We investigate the orbits of $\Gamma$ acting on the two-element subsets of $\mbb{F}_p\cup \mbb{F}_p'$. Note that if $G$ is one of such orbits, then for $r,s\in\mbb{F}_p$, $\{r,s\}\in G$ if and only if $\{r',s'\}\in G$ (because $\alpha$ sends $\{r,s\}$ to $\{r',s'\}$).\\
Let $G_0$ be the orbit $\{0,0'\}$. The orbit of $\{0,0'\}$ under $\beta$ gives us all edges of the form $\{x,0'\}$ with $x\in\mbb{F}_p$. Then, fixing $x\in\mbb{F}_p$, the orbit of $\{x,0'\}$ under $\gamma$ gives all edges of the form $\{x,y'\}$ with $y'\in\mbb{F}_p'$. We conclude that all edges $\{x,y'\}$, $x\in\mbb{F}_p$, $y'\in\mbb{F}_p'$ are in $G_0$. The group $\Gamma$ preserves the set of edges $\{\{x,y'\}: x\in\mbb{F}_p, y'\in\mbb{F}_p'\}$, so $G_0=\{\{x,y'\}: x\in\mbb{F}_p, y'\in\mbb{F}_p'\}$ ($G_0$ is isomorphic to the complete bipartite graph $K_{p,p}$).\\
Let $G$ be any orbit of $\Gamma$ different from $G_0$. Then, $G$ does not contain any edge of the form $\{x,y'\}$ with $x\in\mbb{F}_p$, $y'\in\mbb{F}_p'$. Let $\{r,s\}\in G$, then $r,s\in \mbb{F}_p$ or $r,s\in\mbb{F}_p'$. We can assume that $r,s\in\mbb{F}_p$, since $\{r,s\}\in G$ if and only if $\{r',s'\}\in G$. In order to determine $G$, it suffices to determine the orbit of $\beta$ acting on the two-element subsets of $\mbb{F}_p$. If $K$ represents the orbit of $\{r,s\}$ under $\beta$, then $G=K\cup\alpha(K)$. Now, the permutation $\beta^{p-r}$ sends $\{r,s\}$ to $\{0,s-r\}$. Let $t=s-r\neq 0$. The orbit of $\{0,t\}$ under $\beta$ consists of the $p$ edges $\{0,t\},\{1,t+1\},\{2,t+2\},\ldots,\{p-1,t+p-1\}$, this set of edges is precisely $K$. We want to describe $K$ in a more convenient way. The permutation $\beta^{t}$ is also a generator of $\langle\beta\rangle$ and the orbit of $\{0,t\}$ under $\beta^t$ is describe as the set of edges $\{0,t\},\{t,2t\},\{2t,3t\},\ldots,\{(p-2)t,(p-1)t\},\{(p-1)t,0\}$,  this is $K$. In the notation of lemma \ref{D3}, $K=C(t)$. Thus, we have that $G=G_t=C(t)\cup \alpha(C(t))=C(t)\cup C(t')$ (since $\alpha(C(t))=C(t')$).

There are $(p+1)/2$ orbits of $\Gamma$, $G_0\cong K_{p,p}$ and $C(t)\cup C(t')\cong 2C_p$ for $t=1,2,\ldots,(p-1)/2$.\\
Since $\PP^{\Gamma}\neq\vn$, at least one the orbits of $\Gamma$ belongs to $\PP$ and this ends the proof.
\edem

\bcor \label{dimension} Let $\PP$ be a nontrivial monotone and non-evasive graph property on $2p$ vertices, where $p>3$ is prime. Then, $dim\PP\geq 4p-1$. 
\ecor

\bdem By proposition \ref{2cpokpp}, at least one of $2C_p$, $K_{p,p}$ belongs to $\PP$. If $K_{p,p}$ is in $\PP$, then $\PP$ considered as a simplicial complex contains a face of dimension $p^2-1$. If $\PP$ does not contain $K_{p,p}$, then $\PP$ contains $2C_p$. Therefore, $\PP$ contains all the graphs $C(t)\cup C(t')$ in the proof of lemma \ref{2cpokpp}. Since $(p-1)/2\geq 2$, there are at least two orbits of the form $C(t)\cup C(t')$. By theorem \ref{th2}, $\chi(|\PP^{\Gamma}|=1$, so $\PP$ must contain a graph which is the union of two of the orbits $C(t)\cup C8t')$, this union has $4p$ edges and is a face of $\PP$ of dimension $4p-1$. As $p^2-1\geq 4p-1$, we find that in any case $dim \PP\geq 4p-1$.\edem

\bobs
A result of Bjorner establishes that that for a vertex homogeneous simplicial complex $K$ on a finite set of cardinality $m$ with $\chi(K)=1$, the dimension of $K$ satisfies $dim K\geq M-1$, where $M$ is the maximum prime power dividing $m$ (see \cite{lutz}). Corollary \ref{dimension} says that for a nontrivial monotone and non-evasive property $\PP$ of graphs 0n $2p$ vertices, $dim \PP\geq 4p-1$. $\PP$ represents a simplicial complex on a vertex set of $p(2p-1)$ elements, thus Bjorner's bound for the dimension gives $2p-2$ in the best case and $4p-1$ is a better lower bound.
\eobs

\section{Evasiveness and graphs on Ten Vertices}\label{10vertices}

Now we want to apply similar ideas as in the previous section to nontrivial monotone and non-evasive graph properties $\PP$ on 10 vertices. We will label the ten vertices for our graphs as $0,1,2,3,4,0',1',2',3',4'$. First, we apply lemma \ref{D3} to $p=5$ to find all graphs $G$ on the 5 vertices 0,1,2,3,4, such that 5 does not divide $|[G]|$. They are $\ol{K}_5, C_5\cong C(1)\cong C(2)$ and $K_5=C(1,2)$ (note that $C(2)=\ol{C(1)}$). Correspondingly, for the 5 vertices $0',1',2',3',4'$, we have the graphs $C(1')\cong C(2')$ and $C(1',2')$.\\

Suppose that, for a graph $G$ on 10 vertices, $5$ does not divide $|[G]|$ and apply lemma \ref{D5}. Then, $G$ is isomorphic to $\ol{K}_{10}$, $K_{10}$, or to one of the 10 graphs $G_i, \ol{G}_i$, $i=1,2,3,4,5$ shown in figure \ref{10}.

\begin{figure}
\bcen
\begin{tikzpicture}

\draw (0.25,.75)-- (.75,.5);    
\draw (.75,.5)-- (1,0);            
\draw (1,0)-- (.75,-.5);           
\draw (.75,-.5)-- (0.25,-.75);    
\draw (0.25,.75)-- (0.25,-.75);

\begin{scriptsize}
\draw [fill=black] (0.25,.75) circle (1pt);
\draw (0.257,.88) node {0};
\draw [fill=black] (.75,.5) circle (1pt);
\draw (.9,.517) node {1};
\draw [fill=black] (1,0) circle (1pt);
\draw (1.215,0) node {2};
\draw [fill=black] (.75,-.5) circle (1pt);
\draw (.9,-.517) node {3};
\draw [fill=black] (0.257,-.75) circle (1pt);
\draw (0.257,-.9) node {4};
\draw [fill=black] (-0.25,-.75) circle (1pt);
\draw (-0.257,-.9) node {$4'$};
\draw [fill=black] (-.75,-.5) circle (1pt);
\draw (-.9,-.517) node {$3'$};
\draw [fill=black] (-1,0) circle (1pt);
\draw (-1.215,0) node {$2'$};
\draw [fill=black] (-.75,.5) circle (1pt);
\draw (-.9,.517) node {$1'$};
\draw [fill=black] (-0.25,.75) circle (1pt);
\draw (-0.257,.88) node {$0'$};
\draw (0,-1.4) node {$G_1=\ol{K}_5\cup C(1)$};
\end{scriptsize}
\end{tikzpicture}
\begin{tikzpicture}

\draw (0.25,.75)-- (.75,.5);    
\draw (.75,.5)-- (1,0);            
\draw (1,0)-- (.75,-.5);           
\draw (.75,-.5)-- (0.25,-.75);    
\draw (0.25,.75)-- (0.25,-.75);

\draw (0.25,.75)-- (1,0);        
\draw (0.25,.75)-- (.75,-.5);   
\draw (.75,.5)-- (.75,-.5);       
\draw (.75,.5)-- (0.25,-.75);    
\draw (1,0)-- (0.25,-.75);        

\begin{scriptsize}
\draw [fill=black] (0.25,.75) circle (1pt);
\draw (0.257,.88) node {0};
\draw [fill=black] (.75,.5) circle (1pt);
\draw (.9,.517) node {1};
\draw [fill=black] (1,0) circle (1pt);
\draw (1.215,0) node {2};
\draw [fill=black] (.75,-.5) circle (1pt);
\draw (.9,-.517) node {3};
\draw [fill=black] (0.257,-.75) circle (1pt);
\draw (0.257,-.9) node {4};
\draw [fill=black] (-0.25,-.75) circle (1pt);
\draw (-0.257,-.9) node {$4'$};
\draw [fill=black] (-.75,-.5) circle (1pt);
\draw (-.9,-.517) node {$3'$};
\draw [fill=black] (-1,0) circle (1pt);
\draw (-1.215,0) node {$2'$};
\draw [fill=black] (-.75,.5) circle (1pt);
\draw (-.9,.517) node {$1'$};
\draw [fill=black] (-0.25,.75) circle (1pt);
\draw (-0.257,.88) node {$0'$};
\draw (0,-1.4) node {$G_2=\ol{K}_5\cup C(1,2)$};
\end{scriptsize}
\end{tikzpicture}
\begin{tikzpicture}

\draw (0.25,.75)-- (.75,.5);    
\draw (.75,.5)-- (1,0);            
\draw (1,0)-- (.75,-.5);           
\draw (.75,-.5)-- (0.25,-.75);    
\draw (0.25,.75)-- (0.25,-.75);

\draw (-0.25,-.75)-- (-.75,-.5);   
\draw (-0.25,-.75)-- (-0.25,.75);  
\draw (-.75,-.5)-- (-1,0);           
\draw (-1,0)-- (-.75,.5);            
\draw (-.75,.5)-- (-0.25,.75);         

\begin{scriptsize}
\draw [fill=black] (0.25,.75) circle (1pt);
\draw (0.257,.88) node {0};
\draw [fill=black] (.75,.5) circle (1pt);
\draw (.9,.517) node {1};
\draw [fill=black] (1,0) circle (1pt);
\draw (1.215,0) node {2};
\draw [fill=black] (.75,-.5) circle (1pt);
\draw (.9,-.517) node {3};
\draw [fill=black] (0.257,-.75) circle (1pt);
\draw (0.257,-.9) node {4};
\draw [fill=black] (-0.25,-.75) circle (1pt);
\draw (-0.257,-.9) node {$4'$};
\draw [fill=black] (-.75,-.5) circle (1pt);
\draw (-.9,-.517) node {$3'$};
\draw [fill=black] (-1,0) circle (1pt);
\draw (-1.215,0) node {$2'$};
\draw [fill=black] (-.75,.5) circle (1pt);
\draw (-.9,.517) node {$1'$};
\draw [fill=black] (-0.25,.75) circle (1pt);
\draw (-0.257,.88) node {$0'$};
\draw (0,-1.4) node {$G_3=C(1')\cup C(1)$};
\end{scriptsize}
\end{tikzpicture}
\begin{tikzpicture}

\draw (0.25,.75)-- (.75,.5);    
\draw (.75,.5)-- (1,0);            
\draw (1,0)-- (.75,-.5);           
\draw (.75,-.5)-- (0.25,-.75);    
\draw (0.25,.75)-- (0.25,-.75);

\draw (0.25,.75)-- (1,0);        
\draw (0.25,.75)-- (.75,-.5);   
\draw (.75,.5)-- (.75,-.5);       
\draw (.75,.5)-- (0.25,-.75);    
\draw (1,0)-- (0.25,-.75);        

\draw (-0.25,-.75)-- (-.75,-.5);   
\draw (-0.25,-.75)-- (-0.25,.75);  
\draw (-.75,-.5)-- (-1,0);           
\draw (-1,0)-- (-.75,.5);            
\draw (-.75,.5)-- (-0.25,.75);         

\begin{scriptsize}
\draw [fill=black] (0.25,.75) circle (1pt);
\draw (0.257,.88) node {0};
\draw [fill=black] (.75,.5) circle (1pt);
\draw (.9,.517) node {1};
\draw [fill=black] (1,0) circle (1pt);
\draw (1.215,0) node {2};
\draw [fill=black] (.75,-.5) circle (1pt);
\draw (.9,-.517) node {3};
\draw [fill=black] (0.257,-.75) circle (1pt);
\draw (0.257,-.9) node {4};
\draw [fill=black] (-0.25,-.75) circle (1pt);
\draw (-0.257,-.9) node {$4'$};
\draw [fill=black] (-.75,-.5) circle (1pt);
\draw (-.9,-.517) node {$3'$};
\draw [fill=black] (-1,0) circle (1pt);
\draw (-1.215,0) node {$2'$};
\draw [fill=black] (-.75,.5) circle (1pt);
\draw (-.9,.517) node {$1'$};
\draw [fill=black] (-0.25,.75) circle (1pt);
\draw (-0.257,.88) node {$0'$};
\draw (0,-1.4) node {$G_4=C(1')\cup C(1,2)$};
\end{scriptsize}
\end{tikzpicture}
\ecen
\bcen
\begin{tikzpicture}

\draw (0.25,.75)-- (.75,.5);    
\draw (.75,.5)-- (1,0);            
\draw (1,0)-- (.75,-.5);           
\draw (.75,-.5)-- (0.25,-.75);    
\draw (0.25,.75)-- (0.25,-.75);

\draw (0.25,.75)-- (1,0);        
\draw (0.25,.75)-- (.75,-.5);   
\draw (.75,.5)-- (.75,-.5);       
\draw (.75,.5)-- (0.25,-.75);    
\draw (1,0)-- (0.25,-.75);        

\draw (-0.25,-.75)-- (-.75,-.5);   
\draw (-0.25,-.75)-- (-0.25,.75);  
\draw (-.75,-.5)-- (-1,0);           
\draw (-1,0)-- (-.75,.5);            
\draw (-.75,.5)-- (-0.25,.75);         

\draw (-0.25,-.75)-- (-1,0);        
\draw (-0.25,-.75)-- (-.75,.5);     
\draw (-.75,-.5)-- (-.75,.5);        
\draw (-.75,-.5)-- (-0.25,.75);    
\draw (-1,0)-- (-0.25,.75);         

\begin{scriptsize}
\draw [fill=black] (0.25,.75) circle (1pt);
\draw (0.257,.88) node {0};
\draw [fill=black] (.75,.5) circle (1pt);
\draw (.9,.517) node {1};
\draw [fill=black] (1,0) circle (1pt);
\draw (1.215,0) node {2};
\draw [fill=black] (.75,-.5) circle (1pt);
\draw (.9,-.517) node {3};
\draw [fill=black] (0.257,-.75) circle (1pt);
\draw (0.257,-.9) node {4};
\draw [fill=black] (-0.25,-.75) circle (1pt);
\draw (-0.257,-.9) node {$4'$};
\draw [fill=black] (-.75,-.5) circle (1pt);
\draw (-.9,-.517) node {$3'$};
\draw [fill=black] (-1,0) circle (1pt);
\draw (-1.215,0) node {$2'$};
\draw [fill=black] (-.75,.5) circle (1pt);
\draw (-.9,.517) node {$1'$};
\draw [fill=black] (-0.25,.75) circle (1pt);
\draw (-0.257,.88) node {$0'$};
\draw (0,-1.4) node {$G_5=C(1',2')\cup C(1,2)$};
\end{scriptsize}
\end{tikzpicture}
\begin{tikzpicture}

\draw (0.25,.75)-- (1,0);        
\draw (0.25,.75)-- (.75,-.5);   
\draw (.75,.5)-- (.75,-.5);       
\draw (.75,.5)-- (0.25,-.75);    
\draw (1,0)-- (0.25,-.75);        

\draw (-0.25,-.75)-- (-.75,-.5);   
\draw (-0.25,-.75)-- (-0.25,.75);  
\draw (-.75,-.5)-- (-1,0);           
\draw (-1,0)-- (-.75,.5);            
\draw (-.75,.5)-- (-0.25,.75);         

\draw (-0.25,-.75)-- (-1,0);        
\draw (-0.25,-.75)-- (-.75,.5);     
\draw (-.75,-.5)-- (-.75,.5);        
\draw (-.75,-.5)-- (-0.25,.75);    
\draw (-1,0)-- (-0.25,.75);         

\draw (0.25,.75)-- (-0.25,-.75);
\draw (0.25,.75)-- (-.75,-.5);   
\draw (0.25,.75)-- (-1,0);       
\draw (0.25,.75)-- (-.75,.5);    
\draw (0.25,.75)-- (-0.25,.75);   
\draw (.75,.5)-- (-0.25,-.75);   
\draw (.75,.5)-- (-.75,-.5);      
\draw (.75,.5)-- (-1,0);           
\draw (.75,.5)-- (-.75,.5);       
\draw (.75,.5)-- (-0.25,.75);    
\draw (1,0)-- (-0.25,-.75);       
\draw (1,0)-- (-.75,-.5);          
\draw (1,0)-- (-1,0);              
\draw (1,0)-- (-.75,.5);           
\draw (1,0)-- (-0.25,.75);        
\draw (.75,-.5)-- (-0.25,-.75);   
\draw (.75,-.5)-- (-.75,-.5);      
\draw (.75,-.5)-- (-1,0);           
\draw (.75,-.5)-- (-.75,.5);        
\draw (.75,-.5)-- (-0.25,.75);     
\draw (0.25,-.75)-- (-0.25,-.75); 
\draw (0.25,-.75)-- (-.75,-.5);    
\draw (0.25,-.75)-- (-1,0);         
\draw (0.25,-.75)-- (-.75,.5);      
\draw (0.25,-.75)-- (-0.25,.75);   

\begin{scriptsize}
\draw [fill=black] (0.25,.75) circle (1pt);
\draw (0.257,.88) node {0};
\draw [fill=black] (.75,.5) circle (1pt);
\draw (.9,.517) node {1};
\draw [fill=black] (1,0) circle (1pt);
\draw (1.215,0) node {2};
\draw [fill=black] (.75,-.5) circle (1pt);
\draw (.9,-.517) node {3};
\draw [fill=black] (0.257,-.75) circle (1pt);
\draw (0.257,-.9) node {4};
\draw [fill=black] (-0.25,-.75) circle (1pt);
\draw (-0.257,-.9) node {$4'$};
\draw [fill=black] (-.75,-.5) circle (1pt);
\draw (-.9,-.517) node {$3'$};
\draw [fill=black] (-1,0) circle (1pt);
\draw (-1.215,0) node {$2'$};
\draw [fill=black] (-.75,.5) circle (1pt);
\draw (-.9,.517) node {$1'$};
\draw [fill=black] (-0.25,.75) circle (1pt);
\draw (-0.257,.88) node {$0'$};
\draw (0,-1.4) node {$\ol{G}_1=C(1',2')+C(2)$};
\end{scriptsize}
\end{tikzpicture}
\begin{tikzpicture}

\draw (-0.25,-.75)-- (-.75,-.5);   
\draw (-0.25,-.75)-- (-0.25,.75);  
\draw (-.75,-.5)-- (-1,0);           
\draw (-1,0)-- (-.75,.5);            
\draw (-.75,.5)-- (-0.25,.75);         

\draw (-0.25,-.75)-- (-1,0);        
\draw (-0.25,-.75)-- (-.75,.5);     
\draw (-.75,-.5)-- (-.75,.5);        
\draw (-.75,-.5)-- (-0.25,.75);    
\draw (-1,0)-- (-0.25,.75);         

\draw (0.25,.75)-- (-0.25,-.75);
\draw (0.25,.75)-- (-.75,-.5);   
\draw (0.25,.75)-- (-1,0);       
\draw (0.25,.75)-- (-.75,.5);    
\draw (0.25,.75)-- (-0.25,.75);   
\draw (.75,.5)-- (-0.25,-.75);   
\draw (.75,.5)-- (-.75,-.5);      
\draw (.75,.5)-- (-1,0);           
\draw (.75,.5)-- (-.75,.5);       
\draw (.75,.5)-- (-0.25,.75);    
\draw (1,0)-- (-0.25,-.75);       
\draw (1,0)-- (-.75,-.5);          
\draw (1,0)-- (-1,0);              
\draw (1,0)-- (-.75,.5);           
\draw (1,0)-- (-0.25,.75);        
\draw (.75,-.5)-- (-0.25,-.75);   
\draw (.75,-.5)-- (-.75,-.5);      
\draw (.75,-.5)-- (-1,0);           
\draw (.75,-.5)-- (-.75,.5);        
\draw (.75,-.5)-- (-0.25,.75);     
\draw (0.25,-.75)-- (-0.25,-.75); 
\draw (0.25,-.75)-- (-.75,-.5);    
\draw (0.25,-.75)-- (-1,0);         
\draw (0.25,-.75)-- (-.75,.5);      
\draw (0.25,-.75)-- (-0.25,.75);   

\begin{scriptsize}
\draw [fill=black] (0.25,.75) circle (1pt);
\draw (0.257,.88) node {0};
\draw [fill=black] (.75,.5) circle (1pt);
\draw (.9,.517) node {1};
\draw [fill=black] (1,0) circle (1pt);
\draw (1.215,0) node {2};
\draw [fill=black] (.75,-.5) circle (1pt);
\draw (.9,-.517) node {3};
\draw [fill=black] (0.257,-.75) circle (1pt);
\draw (0.257,-.9) node {4};
\draw [fill=black] (-0.25,-.75) circle (1pt);
\draw (-0.257,-.9) node {$4'$};
\draw [fill=black] (-.75,-.5) circle (1pt);
\draw (-.9,-.517) node {$3'$};
\draw [fill=black] (-1,0) circle (1pt);
\draw (-1.215,0) node {$2'$};
\draw [fill=black] (-.75,.5) circle (1pt);
\draw (-.9,.517) node {$1'$};
\draw [fill=black] (-0.25,.75) circle (1pt);
\draw (-0.257,.88) node {$0'$};
\draw (0,-1.4) node {$\ol{G}_2=C(1',2')+\ol{K}_5$};
\end{scriptsize}
\end{tikzpicture}
\begin{tikzpicture}

\draw (0.25,.75)-- (1,0);        
\draw (0.25,.75)-- (.75,-.5);   
\draw (.75,.5)-- (.75,-.5);       
\draw (.75,.5)-- (0.25,-.75);    
\draw (1,0)-- (0.25,-.75);        

\draw (-0.25,-.75)-- (-1,0);        
\draw (-0.25,-.75)-- (-.75,.5);     
\draw (-.75,-.5)-- (-.75,.5);        
\draw (-.75,-.5)-- (-0.25,.75);    
\draw (-1,0)-- (-0.25,.75);         

\draw (0.25,.75)-- (-0.25,-.75);
\draw (0.25,.75)-- (-.75,-.5);   
\draw (0.25,.75)-- (-1,0);       
\draw (0.25,.75)-- (-.75,.5);    
\draw (0.25,.75)-- (-0.25,.75);   
\draw (.75,.5)-- (-0.25,-.75);   
\draw (.75,.5)-- (-.75,-.5);      
\draw (.75,.5)-- (-1,0);           
\draw (.75,.5)-- (-.75,.5);       
\draw (.75,.5)-- (-0.25,.75);    
\draw (1,0)-- (-0.25,-.75);       
\draw (1,0)-- (-.75,-.5);          
\draw (1,0)-- (-1,0);              
\draw (1,0)-- (-.75,.5);           
\draw (1,0)-- (-0.25,.75);        
\draw (.75,-.5)-- (-0.25,-.75);   
\draw (.75,-.5)-- (-.75,-.5);      
\draw (.75,-.5)-- (-1,0);           
\draw (.75,-.5)-- (-.75,.5);        
\draw (.75,-.5)-- (-0.25,.75);     
\draw (0.25,-.75)-- (-0.25,-.75); 
\draw (0.25,-.75)-- (-.75,-.5);    
\draw (0.25,-.75)-- (-1,0);         
\draw (0.25,-.75)-- (-.75,.5);      
\draw (0.25,-.75)-- (-0.25,.75);   

\begin{scriptsize}
\draw [fill=black] (0.25,.75) circle (1pt);
\draw (0.257,.88) node {0};
\draw [fill=black] (.75,.5) circle (1pt);
\draw (.9,.517) node {1};
\draw [fill=black] (1,0) circle (1pt);
\draw (1.215,0) node {2};
\draw [fill=black] (.75,-.5) circle (1pt);
\draw (.9,-.517) node {3};
\draw [fill=black] (0.257,-.75) circle (1pt);
\draw (0.257,-.9) node {4};
\draw [fill=black] (-0.25,-.75) circle (1pt);
\draw (-0.257,-.9) node {$4'$};
\draw [fill=black] (-.75,-.5) circle (1pt);
\draw (-.9,-.517) node {$3'$};
\draw [fill=black] (-1,0) circle (1pt);
\draw (-1.215,0) node {$2'$};
\draw [fill=black] (-.75,.5) circle (1pt);
\draw (-.9,.517) node {$1'$};
\draw [fill=black] (-0.25,.75) circle (1pt);
\draw (-0.257,.88) node {$0'$};
\draw (0,-1.4) node {$\ol{G}_3=C(2')+C(2)$};
\end{scriptsize}
\end{tikzpicture}
\ecen
\bcen
\begin{tikzpicture}

\draw (-0.25,-.75)-- (-1,0);        
\draw (-0.25,-.75)-- (-.75,.5);     
\draw (-.75,-.5)-- (-.75,.5);        
\draw (-.75,-.5)-- (-0.25,.75);    
\draw (-1,0)-- (-0.25,.75);         

\draw (0.25,.75)-- (-0.25,-.75);
\draw (0.25,.75)-- (-.75,-.5);   
\draw (0.25,.75)-- (-1,0);       
\draw (0.25,.75)-- (-.75,.5);    
\draw (0.25,.75)-- (-0.25,.75);   
\draw (.75,.5)-- (-0.25,-.75);   
\draw (.75,.5)-- (-.75,-.5);      
\draw (.75,.5)-- (-1,0);           
\draw (.75,.5)-- (-.75,.5);       
\draw (.75,.5)-- (-0.25,.75);    
\draw (1,0)-- (-0.25,-.75);       
\draw (1,0)-- (-.75,-.5);          
\draw (1,0)-- (-1,0);              
\draw (1,0)-- (-.75,.5);           
\draw (1,0)-- (-0.25,.75);        
\draw (.75,-.5)-- (-0.25,-.75);   
\draw (.75,-.5)-- (-.75,-.5);      
\draw (.75,-.5)-- (-1,0);           
\draw (.75,-.5)-- (-.75,.5);        
\draw (.75,-.5)-- (-0.25,.75);     
\draw (0.25,-.75)-- (-0.25,-.75); 
\draw (0.25,-.75)-- (-.75,-.5);    
\draw (0.25,-.75)-- (-1,0);         
\draw (0.25,-.75)-- (-.75,.5);      
\draw (0.25,-.75)-- (-0.25,.75);   

\begin{scriptsize}
\draw [fill=black] (0.25,.75) circle (1pt);
\draw (0.257,.88) node {0};
\draw [fill=black] (.75,.5) circle (1pt);
\draw (.9,.517) node {1};
\draw [fill=black] (1,0) circle (1pt);
\draw (1.215,0) node {2};
\draw [fill=black] (.75,-.5) circle (1pt);
\draw (.9,-.517) node {3};
\draw [fill=black] (0.257,-.75) circle (1pt);
\draw (0.257,-.9) node {4};
\draw [fill=black] (-0.25,-.75) circle (1pt);
\draw (-0.257,-.9) node {$4'$};
\draw [fill=black] (-.75,-.5) circle (1pt);
\draw (-.9,-.517) node {$3'$};
\draw [fill=black] (-1,0) circle (1pt);
\draw (-1.215,0) node {$2'$};
\draw [fill=black] (-.75,.5) circle (1pt);
\draw (-.9,.517) node {$1'$};
\draw [fill=black] (-0.25,.75) circle (1pt);
\draw (-0.257,.88) node {$0'$};
\draw (0,-1.4) node {$\ol{G}_4=C(2')+\ol{K}_5$};
\end{scriptsize}
\end{tikzpicture}
\begin{tikzpicture}

\draw (0.25,.75)-- (-0.25,-.75);
\draw (0.25,.75)-- (-.75,-.5);   
\draw (0.25,.75)-- (-1,0);       
\draw (0.25,.75)-- (-.75,.5);    
\draw (0.25,.75)-- (-0.25,.75);   
\draw (.75,.5)-- (-0.25,-.75);   
\draw (.75,.5)-- (-.75,-.5);      
\draw (.75,.5)-- (-1,0);           
\draw (.75,.5)-- (-.75,.5);       
\draw (.75,.5)-- (-0.25,.75);    
\draw (1,0)-- (-0.25,-.75);       
\draw (1,0)-- (-.75,-.5);          
\draw (1,0)-- (-1,0);              
\draw (1,0)-- (-.75,.5);           
\draw (1,0)-- (-0.25,.75);        
\draw (.75,-.5)-- (-0.25,-.75);   
\draw (.75,-.5)-- (-.75,-.5);      
\draw (.75,-.5)-- (-1,0);           
\draw (.75,-.5)-- (-.75,.5);        
\draw (.75,-.5)-- (-0.25,.75);     
\draw (0.25,-.75)-- (-0.25,-.75); 
\draw (0.25,-.75)-- (-.75,-.5);    
\draw (0.25,-.75)-- (-1,0);         
\draw (0.25,-.75)-- (-.75,.5);      
\draw (0.25,-.75)-- (-0.25,.75);   

\begin{scriptsize}
\draw [fill=black] (0.25,.75) circle (1pt);
\draw (0.257,.88) node {0};
\draw [fill=black] (.75,.5) circle (1pt);
\draw (.9,.517) node {1};
\draw [fill=black] (1,0) circle (1pt);
\draw (1.215,0) node {2};
\draw [fill=black] (.75,-.5) circle (1pt);
\draw (.9,-.517) node {3};
\draw [fill=black] (0.257,-.75) circle (1pt);
\draw (0.257,-.9) node {4};
\draw [fill=black] (-0.25,-.75) circle (1pt);
\draw (-0.257,-.9) node {$4'$};
\draw [fill=black] (-.75,-.5) circle (1pt);
\draw (-.9,-.517) node {$3'$};
\draw [fill=black] (-1,0) circle (1pt);
\draw (-1.215,0) node {$2'$};
\draw [fill=black] (-.75,.5) circle (1pt);
\draw (-.9,.517) node {$1'$};
\draw [fill=black] (-0.25,.75) circle (1pt);
\draw (-0.257,.88) node {$0'$};
\draw (0,-1.4) node {$\ol{G}_5=K_{5,5}$};
\end{scriptsize}
\end{tikzpicture}
\ecen
\caption{}
\label{10}
\end{figure}

If $\PP$ does not contain any of the 10 graphs described by lemma \ref{D5}, then 5 divides $\chi(\PP)$ and so $\PP$ is evasive. The graph $\ol{K}_{10}$ represents the empty simplex so it does not contribute to $\chi(\PP)$. If $\PP$ in nontrivial, then $\PP$ does not contain $K_{10}$. As a consequence $\PP$ must contain some of the graphs in figure \ref{10}. The following table contains the automorphism groups of $G_i$ for $i=1,2,3,4,5$, also the sizes $|[G_i]|$ and the contribution modulo 5 to $\chi(\PP)$ when $G_i\in\PP$, $(-1)^{m_{G_i}-1}|[G_i]|(mod\ 5)$:\\

\bcen
\begin{tabular}{|c|c|c|c|c|}
\hline
$G_i$ & $Aut(G_i)$ & $|[G_i]|$ & $(-1)^{m_{G_i}-1}|[G_i]|(mod\ 5)$ \\
\hline
$G_1$ & $D_{10}\times S_5$ & $2^4\cdot 3^3\cdot 7$&+4\\
\hline
$G_2$ & $S_5\times S_5$ & $2^{2} \cdot 3^{2} \cdot 7$&-2\\
\hline
$G_3$ & $D_5\wr S_2$ & $2^{5} \cdot 3^{4} \cdot 7$&-4\\
\hline
$G_4$ & $S_5\times D_{10}$ & $2^4\cdot 3^3\cdot 7$&+4\\
\hline
$G_5$ & $S_5\wr S_2$ & $2 \cdot 3^{2} \cdot 7$&-1\\
\hline
\end{tabular}
\ecen

Note that if $G$ is a graph on 10 vertices having $m_G$ edges, then $\ol{G}$ has $m_{\ol{G}}=45-m_G$ edges. The Hasse diagram of the isomorphism classes $[G_i]$'s and $[\ol{G_j}]$'s is the following: 

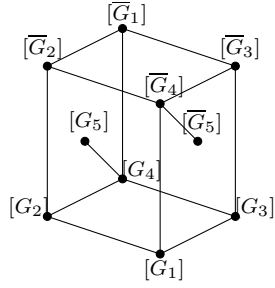
\bcen
\begin{tikzpicture}
\draw (-0.5,2.0)-- (1.0,1.5);
\draw (-0.5,2.0)-- (-1.5,1.5);
\draw (-1.5,1.5)-- (0.0,1.0);
\draw (1.0,1.5)-- (0.0,1.0);
\draw (-0.5,2.0)-- (-0.5,0.0);
\draw (-1.5,1.5)-- (-1.5,-0.5);
\draw (1.0,1.5)-- (1.0,-0.5);
\draw (0.0,1.0)-- (0.0,-1.0);
\draw (-0.5,0.0)-- (-1.5,-0.5);
\draw (-0.5,0.0)-- (1.0,-0.5);
\draw (1.0,-0.5)-- (0.0,-1.0);
\draw (-1.5,-0.5)-- (0.0,-1.0);
\draw (0.0,1.0)-- (0.5,0.5);
\draw (-0.5,0.0)-- (-1.0,0.5);
\begin{scriptsize}
\draw [fill=black] (-0.5,2.0) circle (1.5pt);
\draw (-0.4470185562990762,2.2) node {$[\ol{G}_1]$};
\draw [fill=black] (1.0,1.5) circle (1.5pt);
\draw (1.0594123154026696,1.75) node {$[\ol{G}_3]$};
\draw [fill=black] (-1.5,1.5) circle (1.5pt);
\draw (-1.56020081743392,1.75) node {$[\ol{G}_2]$};
\draw [fill=black] (0.0,1.0) circle (1.5pt);
\draw (0.05248220642308159,1.251101700805588242) node {$[\ol{G}_4]$};
\draw [fill=black] (-0.5,0.0) circle (1.5pt);
\draw (-0.24470185562990762,0.1111685551145078) node {$[G_4]$};
\draw [fill=black] (-1.5,-0.5) circle (1.5pt);
\draw (-1.746020081743392,-0.38833220760765036) node {$[G_2]$};
\draw [fill=black] (1.0,-0.5) circle (1.5pt);
\draw (1.2594123154026696,-0.38833220760765036) node {$[G_3]$};
\draw [fill=black] (0.0,-1.0) circle (1.5pt);
\draw (0.05248220642308159,-1.208878329703298085) node {$[G_1]$};
\draw [fill=black] (0.5,0.5) circle (1.5pt);
\draw (0.5519829691452394,0.75610669317836666) node {$[\ol{G}_5]$};
\draw [fill=black] (-1.0,0.5) circle (1.5pt);
\draw (-0.9465193190212341,0.75610669317836666) node {$[G_5]$};
\end{scriptsize}
\end{tikzpicture}
\ecen

If $\PP$ is non-evasive, then $\PP$ contains some of the isomorphism classes $[G_i], [\ol{G}_j]$. The set of such isomorphism classes contained in $\PP$ becomes an \tit{order ideal} of the poset above. Therefore, $\chi(\PP)$ is congruent modulo $5$ to the sum of the terms $(-1)^{m_G-1}|[G]|$ where $[G]$ belongs to such order ideal. We need to determine all the order ideals of the the poset above and also determine which of them can make $\chi(\PP)$ congruent to 1 $mod\ 5$. There are exactly 9 of these order ideals:
$$\bmat
I_1=\{[\ol{G_5}]\},\\
I_2=\{[G_1],[\ol{G_4}], [\ol{G_5}]\},\\
I_3=\{[G_1],[G_3],[\ol{G_5}]\},\\
I_4=\{[G_1],[G_2], [\ol{G_2}],[\ol{G_4}],[\ol{G_5}]\},\\
I_5=\{[G_1],[G_2],[G_3],[G_4],[G_5]\},\\
I_6=\{[G_1],[G_3],[\ol{G_3}],[\ol{G_4}],[\ol{G_5}]\},\\
I_7=\{[G_1],[G_2],[G_3],[G_4],[\ol{G_2}],[\ol{G_4}],[\ol{G_5}]\},\\
I_8=\{[G_1],[G_2],[G_3],[\ol{G_2}],[\ol{G_3}],[\ol{G_4}],[\ol{G_5}]\},\\
I_9=\{[G_1],[G_2],[G_3],[G_4],[\ol{G_1}],[\ol{G_2}],[\ol{G_3}],[\ol{G_4}],[\ol{G_5}]\}.
\emat$$
 If $\PP$ just contain the isomorphism classes in $I_k$, $k=1,\ldots, 9$, we say that $\PP$ is of \tit{type} $k$. It can be verified that $\PP$ is of type $k$ if and only if $\PP^{\ast}$ is of type $10-k$.\\

Now, we are going to show that types 1, 3, 7 and 9 cannot happen.\\

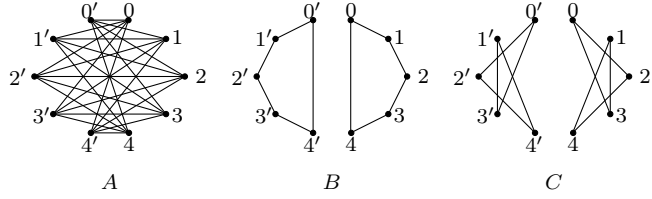
\begin{figure}
\bcen
\begin{tikzpicture}
\draw (0.25,.75)-- (-0.25,-.75);
\draw (0.25,.75)-- (-.75,-.5);   
\draw (0.25,.75)-- (-1,0);       
\draw (0.25,.75)-- (-.75,.5);    
\draw (0.25,.75)-- (-0.25,.75);   
\draw (.75,.5)-- (-0.25,-.75);   
\draw (.75,.5)-- (-.75,-.5);      
\draw (.75,.5)-- (-1,0);           
\draw (.75,.5)-- (-.75,.5);       
\draw (.75,.5)-- (-0.25,.75);    
\draw (1,0)-- (-0.25,-.75);       
\draw (1,0)-- (-.75,-.5);          
\draw (1,0)-- (-1,0);              
\draw (1,0)-- (-.75,.5);           
\draw (1,0)-- (-0.25,.75);        
\draw (.75,-.5)-- (-0.25,-.75);   
\draw (.75,-.5)-- (-.75,-.5);      
\draw (.75,-.5)-- (-1,0);           
\draw (.75,-.5)-- (-.75,.5);        
\draw (.75,-.5)-- (-0.25,.75);     
\draw (0.25,-.75)-- (-0.25,-.75); 
\draw (0.25,-.75)-- (-.75,-.5);    
\draw (0.25,-.75)-- (-1,0);         
\draw (0.25,-.75)-- (-.75,.5);      
\draw (0.25,-.75)-- (-0.25,.75);   

\begin{scriptsize}
\draw [fill=black] (0.25,.75) circle (1pt);
\draw (0.257,.88) node {0};
\draw [fill=black] (.75,.5) circle (1pt);
\draw (.9,.517) node {1};
\draw [fill=black] (1,0) circle (1pt);
\draw (1.215,0) node {2};
\draw [fill=black] (.75,-.5) circle (1pt);
\draw (.9,-.517) node {3};
\draw [fill=black] (0.257,-.75) circle (1pt);
\draw (0.257,-.9) node {4};
\draw [fill=black] (-0.25,-.75) circle (1pt);
\draw (-0.257,-.9) node {$4'$};
\draw [fill=black] (-.75,-.5) circle (1pt);
\draw (-.9,-.517) node {$3'$};
\draw [fill=black] (-1,0) circle (1pt);
\draw (-1.215,0) node {$2'$};
\draw [fill=black] (-.75,.5) circle (1pt);
\draw (-.9,.517) node {$1'$};
\draw [fill=black] (-0.25,.75) circle (1pt);
\draw (-0.257,.88) node {$0'$};
\draw (0,-1.4) node {$A$};
\end{scriptsize}
\end{tikzpicture}
\begin{tikzpicture}
\draw (0.25,.75)-- (.75,.5);    
\draw (.75,.5)-- (1,0);            
\draw (1,0)-- (.75,-.5);           
\draw (.75,-.5)-- (0.25,-.75);    
\draw (0.25,.75)-- (0.25,-.75);
\draw (-0.25,-.75)-- (-.75,-.5);   
\draw (-0.25,-.75)-- (-0.25,.75);  
\draw (-.75,-.5)-- (-1,0);           
\draw (-1,0)-- (-.75,.5);            
\draw (-.75,.5)-- (-0.25,.75);         

\begin{scriptsize}
\draw [fill=black] (0.25,.75) circle (1pt);
\draw (0.257,.88) node {0};
\draw [fill=black] (.75,.5) circle (1pt);
\draw (.9,.517) node {1};
\draw [fill=black] (1,0) circle (1pt);
\draw (1.215,0) node {2};
\draw [fill=black] (.75,-.5) circle (1pt);
\draw (.9,-.517) node {3};
\draw [fill=black] (0.257,-.75) circle (1pt);
\draw (0.257,-.9) node {4};
\draw [fill=black] (-0.25,-.75) circle (1pt);
\draw (-0.257,-.9) node {$4'$};
\draw [fill=black] (-.75,-.5) circle (1pt);
\draw (-.9,-.517) node {$3'$};
\draw [fill=black] (-1,0) circle (1pt);
\draw (-1.215,0) node {$2'$};
\draw [fill=black] (-.75,.5) circle (1pt);
\draw (-.9,.517) node {$1'$};
\draw [fill=black] (-0.25,.75) circle (1pt);
\draw (-0.257,.88) node {$0'$};
\draw (0,-1.4) node {$B$};
\end{scriptsize}
\end{tikzpicture}
\begin{tikzpicture}
\draw (0.25,.75)-- (1,0);        
\draw (0.25,.75)-- (.75,-.5);   
\draw (.75,.5)-- (.75,-.5);       
\draw (.75,.5)-- (0.25,-.75);    
\draw (1,0)-- (0.25,-.75);        
\draw (-0.25,-.75)-- (-1,0);        
\draw (-0.25,-.75)-- (-.75,.5);     
\draw (-.75,-.5)-- (-.75,.5);        
\draw (-.75,-.5)-- (-0.25,.75);    
\draw (-1,0)-- (-0.25,.75);         

\begin{scriptsize}
\draw [fill=black] (0.25,.75) circle (1pt);
\draw (0.257,.88) node {0};
\draw [fill=black] (.75,.5) circle (1pt);
\draw (.9,.517) node {1};
\draw [fill=black] (1,0) circle (1pt);
\draw (1.215,0) node {2};
\draw [fill=black] (.75,-.5) circle (1pt);
\draw (.9,-.517) node {3};
\draw [fill=black] (0.257,-.75) circle (1pt);
\draw (0.257,-.9) node {4};
\draw [fill=black] (-0.25,-.75) circle (1pt);
\draw (-0.257,-.9) node {$4'$};
\draw [fill=black] (-.75,-.5) circle (1pt);
\draw (-.9,-.517) node {$3'$};
\draw [fill=black] (-1,0) circle (1pt);
\draw (-1.215,0) node {$2'$};
\draw [fill=black] (-.75,.5) circle (1pt);
\draw (-.9,.517) node {$1'$};
\draw [fill=black] (-0.25,.75) circle (1pt);
\draw (-0.257,.88) node {$0'$};
\draw (0,-1.4) node {$C$};
\end{scriptsize}
\end{tikzpicture}

\ecen
\caption{Orbits of $\langle(00')(11')(22')(33')(44'),(01234),(0'1'2'3'4')\rangle$}
\label{F2}
\end{figure}

First we show that $\PP$ cannot be of type 3 nor 7. Proposition \ref{2cpokpp} implies that $\PP$ contains one of the 2, $K_{5,5}\cong \ol{G}_5$ or $2C_5\cong G_3$. The corresponding Oliver group is $\Gamma=\langle(00')(11')(22')(33')(44'),(01234),(0'1'2'3'4')\rangle$. The orbits of $\Gamma$, $A\cong \ol{G}_5, B\cong G_3, C\cong G_3$ are shown in figure \ref{F2}. By theorem \ref{th2}, $\chi(|\PP|^{\Gamma})=1$.

If $\PP$ is of type 3 or 7, then $G_3, \ol{G}_5\in\PP$. Thus, we have that $A,B,C\in\PP$ and since $\chi(\PP^{\Gamma})=1$, the graphs $A\cup B\cong A\cup C$ are in $\PP$, but $A\cup B\cong \ol{G}_3$, but $\ol{G}_3$ does not belong to $\PP$ (see the order ideals $I_3$ and $I_7$ above).\\

To show that $\PP$ cannot be of type 1 nor 9, we need the following  result of P. A. Smith \cite{smith}.

\bteo\label{smith} If $\Gamma$ is a $p$-group acting on a $\integer/p$-acyclic complex $K$, then $|K|^{\Gamma}$ is also $\integer/p$-acyclic. 
\eteo

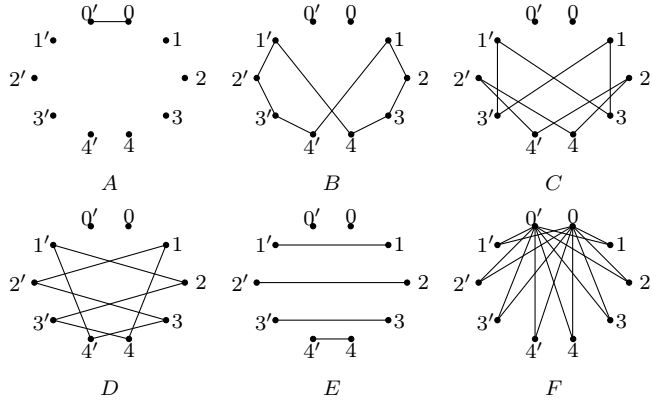
\begin{figure}
\bcen
\begin{tikzpicture}
\draw (0.25,.75)-- (-0.25,.75);   

\begin{scriptsize}
\draw [fill=black] (0.25,.75) circle (1pt);
\draw (0.257,.88) node {0};
\draw [fill=black] (.75,.5) circle (1pt);
\draw (.9,.517) node {1};
\draw [fill=black] (1,0) circle (1pt);
\draw (1.215,0) node {2};
\draw [fill=black] (.75,-.5) circle (1pt);
\draw (.9,-.517) node {3};
\draw [fill=black] (0.257,-.75) circle (1pt);
\draw (0.257,-.9) node {4};
\draw [fill=black] (-0.25,-.75) circle (1pt);
\draw (-0.257,-.9) node {$4'$};
\draw [fill=black] (-.75,-.5) circle (1pt);
\draw (-.9,-.517) node {$3'$};
\draw [fill=black] (-1,0) circle (1pt);
\draw (-1.215,0) node {$2'$};
\draw [fill=black] (-.75,.5) circle (1pt);
\draw (-.9,.517) node {$1'$};
\draw [fill=black] (-0.25,.75) circle (1pt);
\draw (-0.257,.88) node {$0'$};
\draw (0,-1.4) node {$A$};
\end{scriptsize}
\end{tikzpicture}
\begin{tikzpicture}
\draw (.75,.5)-- (1,0);            
\draw (.75,.5)-- (-0.25,-.75);   
\draw (1,0)-- (.75,-.5);           
\draw (.75,-.5)-- (0.25,-.75);    
\draw (0.25,-.75)-- (-.75,.5);      
\draw (-0.25,-.75)-- (-.75,-.5);   
\draw (-.75,-.5)-- (-1,0);           
\draw (-1,0)-- (-.75,.5);            

\begin{scriptsize}
\draw [fill=black] (0.25,.75) circle (1pt);
\draw (0.257,.88) node {0};
\draw [fill=black] (.75,.5) circle (1pt);
\draw (.9,.517) node {1};
\draw [fill=black] (1,0) circle (1pt);
\draw (1.215,0) node {2};
\draw [fill=black] (.75,-.5) circle (1pt);
\draw (.9,-.517) node {3};
\draw [fill=black] (0.257,-.75) circle (1pt);
\draw (0.257,-.9) node {4};
\draw [fill=black] (-0.25,-.75) circle (1pt);
\draw (-0.257,-.9) node {$4'$};
\draw [fill=black] (-.75,-.5) circle (1pt);
\draw (-.9,-.517) node {$3'$};
\draw [fill=black] (-1,0) circle (1pt);
\draw (-1.215,0) node {$2'$};
\draw [fill=black] (-.75,.5) circle (1pt);
\draw (-.9,.517) node {$1'$};
\draw [fill=black] (-0.25,.75) circle (1pt);
\draw (-0.257,.88) node {$0'$};
\draw (0,-1.4) node {$B$};
\end{scriptsize}
\end{tikzpicture}
\begin{tikzpicture}
\draw (.75,.5)-- (.75,-.5);       
\draw (.75,.5)-- (-.75,-.5);      
\draw (1,0)-- (0.25,-.75);        
\draw (1,0)-- (-0.25,-.75);       
\draw (.75,-.5)-- (-.75,.5);        
\draw (0.25,-.75)-- (-1,0);         
\draw (-0.25,-.75)-- (-1,0);        
\draw (-.75,-.5)-- (-.75,.5);        

\begin{scriptsize}
\draw [fill=black] (0.25,.75) circle (1pt);
\draw (0.257,.88) node {0};
\draw [fill=black] (.75,.5) circle (1pt);
\draw (.9,.517) node {1};
\draw [fill=black] (1,0) circle (1pt);
\draw (1.215,0) node {2};
\draw [fill=black] (.75,-.5) circle (1pt);
\draw (.9,-.517) node {3};
\draw [fill=black] (0.257,-.75) circle (1pt);
\draw (0.257,-.9) node {4};
\draw [fill=black] (-0.25,-.75) circle (1pt);
\draw (-0.257,-.9) node {$4'$};
\draw [fill=black] (-.75,-.5) circle (1pt);
\draw (-.9,-.517) node {$3'$};
\draw [fill=black] (-1,0) circle (1pt);
\draw (-1.215,0) node {$2'$};
\draw [fill=black] (-.75,.5) circle (1pt);
\draw (-.9,.517) node {$1'$};
\draw [fill=black] (-0.25,.75) circle (1pt);
\draw (-0.257,.88) node {$0'$};
\draw (0,-1.4) node {$C$};
\end{scriptsize}
\end{tikzpicture}

\begin{tikzpicture}
\draw (.75,.5)-- (0.25,-.75);    
\draw (.75,.5)-- (-1,0);           
\draw (1,0)-- (-.75,-.5);          
\draw (1,0)-- (-.75,.5);           
\draw (.75,-.5)-- (-0.25,-.75);   
\draw (.75,-.5)-- (-1,0);           
\draw (0.25,-.75)-- (-.75,-.5);    
\draw (-0.25,-.75)-- (-.75,.5);     

\begin{scriptsize}
\draw [fill=black] (0.25,.75) circle (1pt);
\draw (0.257,.88) node {0};
\draw [fill=black] (.75,.5) circle (1pt);
\draw (.9,.517) node {1};
\draw [fill=black] (1,0) circle (1pt);
\draw (1.215,0) node {2};
\draw [fill=black] (.75,-.5) circle (1pt);
\draw (.9,-.517) node {3};
\draw [fill=black] (0.257,-.75) circle (1pt);
\draw (0.257,-.9) node {4};
\draw [fill=black] (-0.25,-.75) circle (1pt);
\draw (-0.257,-.9) node {$4'$};
\draw [fill=black] (-.75,-.5) circle (1pt);
\draw (-.9,-.517) node {$3'$};
\draw [fill=black] (-1,0) circle (1pt);
\draw (-1.215,0) node {$2'$};
\draw [fill=black] (-.75,.5) circle (1pt);
\draw (-.9,.517) node {$1'$};
\draw [fill=black] (-0.25,.75) circle (1pt);
\draw (-0.257,.88) node {$0'$};
\draw (0,-1.4) node {$D$};
\end{scriptsize}
\end{tikzpicture}
\begin{tikzpicture}
\draw (.75,.5)-- (-.75,.5);       
\draw (1,0)-- (-1,0);              
\draw (.75,-.5)-- (-.75,-.5);      
\draw (0.25,-.75)-- (-0.25,-.75); 

\begin{scriptsize}
\draw [fill=black] (0.25,.75) circle (1pt);
\draw (0.257,.88) node {0};
\draw [fill=black] (.75,.5) circle (1pt);
\draw (.9,.517) node {1};
\draw [fill=black] (1,0) circle (1pt);
\draw (1.215,0) node {2};
\draw [fill=black] (.75,-.5) circle (1pt);
\draw (.9,-.517) node {3};
\draw [fill=black] (0.257,-.75) circle (1pt);
\draw (0.257,-.9) node {4};
\draw [fill=black] (-0.25,-.75) circle (1pt);
\draw (-0.257,-.9) node {$4'$};
\draw [fill=black] (-.75,-.5) circle (1pt);
\draw (-.9,-.517) node {$3'$};
\draw [fill=black] (-1,0) circle (1pt);
\draw (-1.215,0) node {$2'$};
\draw [fill=black] (-.75,.5) circle (1pt);
\draw (-.9,.517) node {$1'$};
\draw [fill=black] (-0.25,.75) circle (1pt);
\draw (-0.257,.88) node {$0'$};
\draw (0,-1.4) node {$E$};
\end{scriptsize}
\end{tikzpicture}
\begin{tikzpicture}
\draw (0.25,.75)-- (.75,.5);    
\draw (0.25,.75)-- (1,0);        
\draw (0.25,.75)-- (.75,-.5);   
\draw (0.25,.75)-- (0.25,-.75);
\draw (0.25,.75)-- (-0.25,-.75);
\draw (0.25,.75)-- (-.75,-.5);   
\draw (0.25,.75)-- (-1,0);       
\draw (0.25,.75)-- (-.75,.5);    
\draw (.75,.5)-- (-0.25,.75);    
\draw (1,0)-- (-0.25,.75);        
\draw (.75,-.5)-- (-0.25,.75);     
\draw (0.25,-.75)-- (-0.25,.75);   
\draw (-0.25,-.75)-- (-0.25,.75);  
\draw (-.75,-.5)-- (-0.25,.75);    
\draw (-1,0)-- (-0.25,.75);         
\draw (-.75,.5)-- (-0.25,.75);         

\begin{scriptsize}
\draw [fill=black] (0.25,.75) circle (1pt);
\draw (0.257,.88) node {0};
\draw [fill=black] (.75,.5) circle (1pt);
\draw (.9,.517) node {1};
\draw [fill=black] (1,0) circle (1pt);
\draw (1.215,0) node {2};
\draw [fill=black] (.75,-.5) circle (1pt);
\draw (.9,-.517) node {3};
\draw [fill=black] (0.257,-.75) circle (1pt);
\draw (0.257,-.9) node {4};
\draw [fill=black] (-0.25,-.75) circle (1pt);
\draw (-0.257,-.9) node {$4'$};
\draw [fill=black] (-.75,-.5) circle (1pt);
\draw (-.9,-.517) node {$3'$};
\draw [fill=black] (-1,0) circle (1pt);
\draw (-1.215,0) node {$2'$};
\draw [fill=black] (-.75,.5) circle (1pt);
\draw (-.9,.517) node {$1'$};
\draw [fill=black] (-0.25,.75) circle (1pt);
\draw (-0.257,.88) node {$0'$};
\draw (0,-1.4) node {$F$};
\end{scriptsize}
\end{tikzpicture}

\ecen
\caption{Orbits of $\langle(00'),(12341'2'3'4')\rangle$}
\label{F1}
\end{figure}

Suppose that $\PP$ is of type 1 and let $\Gamma=\langle(00'), (12341'2'3'4')\rangle$. Then, $\Gamma$ is a $2$-group. The (potential) vertices of $\PP^{\Gamma}$ are the orbits of $\Gamma$ acting on the vertices of the simplicial complex $\PP$. These orbits are shown in figure \ref{F1}. 

Since $\PP$ is of type 1, $\PP$ contains $\ol{G}_5$ but does not contain $G_1$. The graphs in the following list are in $\PP$ because each of them is isomorphic to a subgraph of $\ol{G}_5$:
\[
A,B,C,D,E,A\cup B,A\cup C,A\cup D, A\cup E, B\cup D.\]
The graphs in the following list are not in $\PP$ because each of them contains a subgraph isomorphic to $G_1$:
\[
B\cup C, B\cup E, B\cup F, C\cup D, C\cup F, D\cup E, D\cup F, E\cup F.\]
Observe that the graphs $F, A\cup F$ and $C\cup E$ are not in these lists.\\

Now, $A$ is one of the vertices of the simplicial complex $\PP$ and $A$ is also one of the fixed points of $\Gamma$. Therefore, $\Gamma$ acts on $lk_{\PP}(A)$. Moreover $lk_{\PP}(A)$ is a non-evasive complex. The fixed point set of the action of $\Gamma$ on $lk_{\PP}(A)$ is given by $lk_{\PP}(A)^{\Gamma}=lk_{\PP^{\Gamma}}(A)$.\\
From the lists of graphs shown above we see that $B,C,D,E$ are vertices of $lk_{\PP}(A)^{\Gamma}$. Since $lk_{\PP}(A)$ is $\integer/2$-acyclic, $lk_{\PP}(A)^{\Gamma}$ is $\integer/2$-acyclic by theorem \ref{smith}. Then, $lk_{\PP}(A)^{\Gamma}$ is connected.\\

The graph $F$ cannot be a vertex of the simplicial complex $lk_{\PP}(A)^{\Gamma}$ because on the contrary $F$ would be an isolated vertex of $lk_{\PP}(A)^{\Gamma}$. Thus, $lk_{\PP}(A)^{\Gamma}$ has precisely the vertices $B,C,D,E$, and by the lists of graphs shown above, the only other faces of $lk_{\PP}(A)^{\Gamma}$ can be $\{B,D\}$ and $\{C,E\}$. In any case $lk_{\PP}(A)^{\Gamma}$ results to be non-connected. This contradiction proves that $\PP$ cannot be of type 1.\\

If $\PP$ is of type 9, then $\PP^{\ast}$ is of type 1 and the argument above applies to $\PP^{\ast}$ to conclude that $\PP^{\ast}$ cannot be of type 1. Thus, $\PP$ cannot be of type 9.

\bobs
We have not found the appropriate Oliver groups to use in order to discard the 5 remaining possible types for $\PP$. If a nontrivial monotone and non-evasive graph property on 10 vertices $\PP$ exists, then it must contain the isomorphism classes in one of the order ideals $I_k$ for $k=2,4,5,6,8$.
\eobs

\end{document}